\documentclass[11pt]{amsart}
\usepackage{amsmath,amsthm, amscd, amssymb, amsfonts}
\usepackage[all]{xy}
\usepackage{mathrsfs}
\usepackage{enumitem}
\usepackage{hyperref}
\usepackage[ansinew]{inputenc}
\usepackage{graphicx,fancyhdr}

\usepackage[dvips, dvipsnames, usenames]{color}

\newcommand{\Dtriangle}[6]{
\xymatrix@R-12pt{  &    \overset{#2}{\circ} \ar  @{-}[dl]_{#4}\ar  @{-}[dr]^{#5} & \\
\overset{#1}{\circ} \ar  @{-}[rr]^{#6}  &  &\overset{#3}{\circ} }}

\newcommand{\comment}[1]{}
\newcommand{\yd}[1]{{}^{ #1 }_{ #1 }\mathcal{YD}}

\numberwithin{equation}{section}
\newtheorem{theorem}{Theorem}[section]
\newtheorem{lemma}[theorem]{Lemma}
\newtheorem{coro}[theorem]{Corollary}
\newtheorem{conjecture}[theorem]{Conjecture}
\newtheorem{prop}[theorem]{Proposition}

\newtheorem{claim}{Claim}
\newtheorem*{claimo}{Claim}
\newtheorem{algo}{}

\theoremstyle{definition}

\newtheorem{algorithm}[theorem]{Algorithm}
\newtheorem{definition}[theorem]{Definition}
\newtheorem{example}[theorem]{Example}

\theoremstyle{remark}
\newtheorem{remark}[theorem]{Remark}

\newcommand\pf{\begin{proof}}
\newcommand\epf{\end{proof}}

\newcommand{\ztu}{\zeta^{-1}}
\newcommand{\twist}{\simeq_{\text{tw}}}

\newcommand{\Cla}{\operatorname{Cl}_{\operatorname{ab}}}
\newcommand{\GK}{\operatorname{GK-dim}}
\newcommand{\Inn}{\operatorname{Inn}}
\newcommand{\supp}{\operatorname{supp}}

\newcommand{\heis}[2]{\mathtt{H}_{#1}(#2)}
\newcommand{\ut}[2]{\mathtt{UT}_{#1}(#2)}
\newcommand{\Mtt}{\mathtt{M}}

\newcommand{\ba}{ \mathbf{a}}
\newcommand{\bb}{\mathbf{b}}
\newcommand{\br}{\mathbf{r}}
\newcommand{\bs}{\mathbf{s}}

\newcommand{\bi}{\mathbf{i}}
\newcommand{\bI}{\mathbf{I}}

\newcommand{\trid}{\triangleright}

\newcommand{\ku}{\Bbbk}
\newcommand{\kut}{ \Bbbk^{\times}}

\newcommand{\G}{\mathbb G}

\newcommand{\I}{\mathbb I}

\newcommand{\J}{\mathbb J}
\newcommand{\N}{\mathbb N}
\newcommand{\bp}{\mathbf{p}}
\newcommand{\bq}{\mathbf{q}}
\newcommand\Sb{\mathbb S}

\newcommand{\Z}{\mathbb Z}

\newcommand{\cC}{\mathcal{C}}

\newcommand{\D}{\mathbb{D}}

\newcommand{\cH}{\mathcal{H}}

\newcommand{\cJ}{\mathcal{J}}

\newcommand{\Oc}{\mathcal{O}}

\newcommand{\Ss}{{\mathcal S}}

\newcommand{\Uc}{\mathcal{U}}

\newcommand{\orbg}{\mathcal Z}

\newcommand{\Aut}{\operatorname{Aut}}

\newcommand{\id}{\operatorname{id}}
\newcommand{\Ind}{\operatorname{Ind}\nolimits}
\newcommand{\Indec}{\operatorname{Indec}\nolimits}
\newcommand{\Irr}{\operatorname{Irr}\nolimits}

\newcommand{\gr}{\operatorname{gr}}

\newcommand{\Hom}{\operatorname{Hom}}

\newcommand{\Rep}{\operatorname{Rep}}

\newcommand{\toba}{\mathscr{B}}

\newcommand{\ot}{\otimes}

\newcounter{tabla}\stepcounter{tabla}

\begin{document}

\title[On pointed Hopf algebras over nilpotent groups]{On pointed Hopf algebras over nilpotent groups}

\author[Nicol\'as Andruskiewitsch]
{Nicol\'as Andruskiewitsch}

\address{FaMAF-CIEM (CONICET), Universidad Nacional de C\'ordoba,
Medina A\-llen\-de s/n, Ciudad Universitaria, 5000 C\' ordoba, Rep\'
ublica Argentina.} \email{nicolas.andruskiewitsch@unc.edu.ar}

\thanks{The work was partially supported by CONICET,
Secyt (UNC) and the Alexander von Humboldt Foundation
through the Research Group Linkage Programme}

\begin{abstract}
We classify finite-dimensional Nichols algebras over finite nil\-potent groups of odd order in group-theoretical terms.
The main step is to show that conjugacy classes of such finite groups are either abelian or of type C; 
this property also holds for finite conjugacy classes of finitely generated nilpotent groups whose torsion has odd order.
To extend our approach to the setting of finite GK-dimension, we propose a new Conjecture on racks of type C.
We also prove that the bosonization a of Nichols algebra of a Yetter-Drinfeld module over a group  
whose support is an infinite conjugacy class has infinite GK-dimension. 
We apply this to the study of the finite GK-dimensional pointed Hopf algebras over finitely generated
torsion-free nilpotent groups.
\end{abstract}

\maketitle


\section*{Introduction}\label{section:introduction}

\subsection{The context}\label{subsection:intro-context}
Let $\ku$ be an algebraically closed field  of characteristic 0.
This paper contributes to the classification of Hopf algebras with finite Gelfand-Kirillov dimension, $\GK$ for short.
Despite recent interest on this question, see \cite{brown-couto,brown-zhang} and references therein, 
the general structure of such Hopf algebras is still mysterious, 
so it is justified to focus on the class of \emph{pointed} Hopf algebras; 
under this assumption we may follow the method from \cite{AS Pointed HA}, 
first applied in the $\GK$ context in \cite{AS-crelle}, cf. also \cite{R quantum groups}.
By a celebrated theorem of Gromov \cite{gromov}, a finitely generated group has finite growth
if and only if it is nilpotent-by-finite.
Thus, the first major goal within the method of \cite{AS Pointed HA} is to classify Nichols algebras with finite $\GK$ over 
finitely generated nilpotent-by-finite groups. 
Towards this goal, the first natural question is to deal with the classification of 
Nichols algebras with finite $\GK$ over abelian groups.
We collect information on this question needed for the general case.
There are various subclasses to consider:

\begin{itemize}[leftmargin=*,label=$\circ$]
\item Nichols algebras of diagonal type, corresponding to semisimple Yetter-Drinfeld modules over 
finitely generated abelian groups. 
Those with finite dimension were classified in \cite{H-classif}.
Conjecture \ref{conj:AAH} below from \cite{AAH-triang} reduces the classification of the finitely generated  
Nichols algebras of diagonal type with finite $\GK$ to \cite{H-classif}. 
Also \cite{AAH-infinite-rank} deals with Nichols algebras of diagonal type with finite $\GK$ which
are not finitely generated.

\medbreak
\item Next comes the class of Nichols algebras of \emph{blocks \& points}; here the classification of those with finite $\GK$
was achieved in \cite{AAH-triang} for finite, and \cite{AAH-infinite-rank} for infinite, rank;
both assume the validity of Conjecture \ref{conj:AAH}.

\medbreak
\item There are decomposable Yetter-Drinfeld modules over abelian groups that are not of the form 
\emph{blocks \& points}; they contain components known as \emph{pale blocks}. 
The classification of those with finite $\GK$ in rank 3 was also obtained in \cite{AAH-triang}
while rank 4 is work in progress \cite{aam}.
\end{itemize}

Summarizing, in order to have a classification of the 
finitely generated Nichols algebras with finite $\GK$ over abelian groups
 it remains to conclude the classification of the \emph{blocks \& pale blocks \& points} giving rise to Nichols algebras with finite $\GK$
and to prove Conjecture \ref{conj:AAH}; we believe that both objectives could be attained soon.
We should also mention that the classification of finite-dimensional Nichols algebras 
over finite groups is far from complete notwithstanding intense activity in this direction.
See  \cite{A-leyva,Carnovale-Costantini-VI,HV-rank>2} and references therein.

\medbreak
The focus of this article is on Nichols algebras over 
nilpotent groups whose  bosonizations 
have finite $\GK$. The main results of this paper are:

\begin{itemize}[leftmargin=*,label=$\circ$]
\item  The description of the finite-dimensional 
Nichols algebras over a finite nilpotent group of odd order $G$ 
up to the knowledge of the conjugacy classes and the representations of the centralizers of $G$,
Theorem \ref{thm:main-nichols-nilpotent-odd}.

\medbreak
\item  The description
of the  bosonizations $\toba(M) \# \ku G$ with $G$ a torsion-free nilpotent group and $M \in \yd{\ku G}$ semisimple having finite $\GK$ 
up to knowledge of the irreducible representations of $G$. See Theorem \ref{thm:main-nichols-nilpotent-torsionfree}.
\end{itemize}

We next describe  how we achieve these results.

\subsection{Finite nilpotent  groups}
Let us start with finite nilpotent groups of odd order.
Our first basic result, Theorem \ref{thm:odd-order-abelian-or-typeC}, states
that any conjugacy class of such a group is either abelian or of type C. 
The notion of rack of type C was introduced in \cite{ACG-III}, where it was shown
that any Nichols algebra whose support is of type C has infinite dimension.   
Thus we are reduced to deal with abelian conjugacy classes which
give rise to braided vector spaces of diagonal type.
To give a more precise answer, we generalize a technique from \cite{YZ},
the only reference we know on Nichols algebras over nilpotent groups (beyond the abelian case). 
A similar analysis could be carried out in the context of finite $\GK$,
but we need to assume the already mentioned Conjecture \ref{conj:AAH} and 
the new Conjecture \ref{conjecture:typeC-X-finite-GK}
extending the criterium of type C from \cite{ACG-III} to the setting of $\GK$.

\subsection{Finitely generated nilpotent  groups}
Our second basic result, Theorem \ref{thm:infinite-conjugacy-class}, shows that
the bosonization $\toba(M) \# \ku G$ of the Nichols algebra of any Yetter-Drinfeld $M$
module over any finitely generated group $G$, whose support is an infinite conjugacy class, has infinite $\GK$
even if $\toba(M)$ could have finite $\GK$. It is known that any conjugacy class of a 
finitely generated torsion-free nilpotent group $G$ is either infinite or central, 
thus for such groups we just need to study Nichols algebras with central support. 
These arise also over abelian groups, discussed above in \S \ref{subsection:intro-context}. 
For illustration  we list those corresponding to $M$ semisimple,
see Theorem \ref{thm:main-nichols-nilpotent-torsionfree}.

\medbreak
Finally let $G$ be a finitely generated nilpotent group and assume that its torsion subgroup $T$ has odd order.
Then we show that any finite conjugacy class is either abelian or of type C, see Proposition \ref{prop:noncentral-conjugacyclasses-torsion-typeC}, extending Theorem \ref{thm:odd-order-abelian-or-typeC}.
To proceed further we need the validity of Conjectures \ref{conj:AAH} and \ref{conjecture:typeC-X-finite-GK}.
We also make a reduction when the order of $T$ is coprime to 6.

\medbreak The paper has four sections; in the first one we collect preliminary information on Nichols algebras.
Section \ref{section:hopf-nilpotent} deals with conjugacy classes including the basic theorems mentioned above.
In Section \ref{section:nichols-nilpotent} we establish the main results, stating some auxiliary lemmas in more  generality for future applications and discussing a few examples. Comments on the open questions on Hopf algebras over
nilpotent-by-finite groups are in the last Section \ref{section:conclusions}.

\subsection*{Acknowledgment} I thank Caleb Eckhardt and Pavel Shumyatsky
for pointing out proofs of 
Lemmas \ref{lema:noncentral-conjugacyclasses} and \ref{lema:order-conjugacy-class} to me, respectively.
I am also grateful to Sonia Natale for many interesting discussions (on the matters of this paper).

\section{Preliminaries}\label{section:preliminaries}

\subsection{Notations}\label{subsection:notations}
We denote the cardinal of a set $X$ by $\vert X\vert$. 
If $k < \ell$ are non-negative integers, then we set $\I_{k, \ell} = \{i\in \N: k\le i\leq \ell\}$ and just
$\I_{\ell} = \I_{1, \ell}$.
Given a positive integer $\ell$, 
we denote by $\G_\ell$ the group of $\ell$-th roots of unity in $\ku$, and by $\G'_\ell \subset \G_\ell$ the subset of those of order $\ell$. 
The group of all roots of unity is denoted by $\G_{\infty}$ and $\G'_{\infty} := \G_{\infty} - \{1\}$.

Let $G$ be a group. The identity, the group of characters and the center of $G$  are denoted by $e$,   
$\widehat G = \Hom_{\text{groups}}(G, \ku^{\times})$ and $Z(G)$.
The notations $F\leqslant G$, or $G\geqslant F$,  mean that $F$ is a subgroup of $G$, while $F \lhd  G$, or $G \rhd  F$, mean that $F\leqslant G$ is normal. We shall use the notations $x\trid y = xyx^{-1}$,  
$[x, y] = xyx^{-1}y^{-1}$ (the commutator), $\vert x\vert =$ the order of $x$, for $x, y \in G$.
Given $x\in G$, let $\mathcal O_x$ be its conjugacy class and let $G^x$ be its centralizer. 
If emphasis is needed, then we write $\mathcal O_x^G = \mathcal O_x$. 

The symmetric and exterior algebras of a vector space $V$ are denoted $S(V)$ and $\Lambda(V)$ respectively.

Let $\Irr \cC$ be the set of isomorphism classes of simple objects in an abelian category $\cC$.
If $A$ is an algebra and $\cC$ is the category of $A$-modules, then $\Irr A :=\Irr \cC$;
if $A =\ku G$, then $\Irr G :=\Irr A$. Also  $\Indec \cC$ denotes 
the set of isomorphism classes of indecomposable objects in $\cC$ and corespondingly we have 
$\Indec A$, $\Indec G$.

\subsection{Yetter-Drinfeld modules}\label{subsection:yd}
A \emph{braided vector space} is a pair $(V,c)$ where $V$ is a vector space and $c \in GL(V^{\ot 2})$ satisfies 
the braid equation
$(c\ot \id) (\id \ot c) (c\ot \id) = (\id \ot c) (c\ot \id) (\id \ot c)$.
A systematic way of producing braided vector spaces is through 
Yetter-Drinfeld modules over a Hopf algebra $\cH$  (always assumed with bijective antipode);
these are $\cH$-modules and $\cH$-comodules subject to a compatibility condition, see \cite{Rad-libro}.
The category $\yd{\cH}$ of Yetter-Drinfeld modules over $\cH$ is a braided monoidal one.
Hence the notion of  Hopf algebras in $\yd{\cH}$ is available. 
We refer to \cite{A-leyva} for the concepts of Nichols algebra of a Yetter-Drinfeld module (a special kind of Hopf algebra in $\yd{\cH}$) and Nichols algebra of a braided vector space (a non-categorical version of the former), 
central in the approach to the classification of pointed Hopf algebras pursued in this paper.

If $V \in \yd{\cH}$, or if $(V,c)$ is a braided vector space,
then $\toba(V)$ denotes its Nichols algebra  and $\cJ = \cJ(V)$ its ideal of defining relations.

\begin{example}\label{exa:YD-pair}
Yetter-Drinfeld modules of dimension 1 are classified by YD-pairs over $\cH$, that is pairs
$(g, \chi) \in G(H) \times \Hom_{\text{Alg}} (H, \ku)$ satisfying 
\begin{align*}
\chi(h) g &= \chi(h_{(2)}) h_{(1)} g \Ss(h_{(3)}), && h \in H.
\end{align*}
Given a YD-pair $(g, \chi)$, we denote by  $\ku^{\chi}_g \in \yd{\cH}$ the one-dimensional vector space
where $\cH$ acts by $\chi$ and  co-acts  by $g$. Let $q = \chi(g)$. It is well-known that 
$\toba(\ku^{\chi}_g) \simeq \begin{cases}
\ku[T]/T^N, &\text{if } q \in \G'_N, N >1; \\
\ku[T], &\text{ otherwise,}
\end{cases}$ where $T$ is an indeterminate.
\end{example}

\subsection{Nichols algebras of diagonal type}\label{subsection:Nichols}

Let $\theta \in \N$ and $\I := \I_{\theta}$. Given a matrix $\bq = (q_{ij})_{i, j \in \I}\in (\ku^{\times})^{\I \times \I}$,
we denote by $(V, c^{\bq})$ the braided vector space \emph{of diagonal type} 
associated to $\bq$, where $V$ has a basis $(x_i) _ {i \in \I}$ and
\begin{align}\label{eq:def-diagonal}
c^{\bq}(x_i \ot x_j) = q_{ij} x_j \ot x_i, && i, j \in \I. 
\end{align}
In this case we set  $\cJ_{\bq} = \cJ(V)$, $\toba_{\bq} =\toba(V)$, etc.
The \emph{Dynkin diagram} associated to $\bq$ is the  graph with $\theta$ vertices,
where the vertex $i$ is labelled by $q_{ii}$, 
and there is an edge between $i$ and $j$ labelled by $\widetilde{q}_{ij}:= q_{ij}q_{ji}$. 
When $\widetilde{q}_{ij} = 1$, the edge is omitted except sometimes for the needs of the exposition.

\medbreak
Assume that a matrix $\bp = (p_{ij})_{i, j \in \I}$ is \emph{twist-equivalent} to $\bq$, that is
they have the same Dynkin diagram, i.~e. $p_{ii} = q_{ii}$ and $\widetilde{p}_{ij} = \widetilde{q}_{ij}$
for all $i\neq j$. Then the Nichols algebras $\toba_{\bp}$ and $\toba_{\bq}$, which are not necessarily isomorphic, 
have the same Hilbert series, hence the same $\GK$ by \cite[Lemma 6.1]{KL}. 
We shall express this situation by $\toba_{\bp} \twist \toba_{\bq}$.

\medbreak
We may refer to  the connected components of the Dynkin diagram and a fortiori of $\bq$.
For many purposes we may assume that $\bq$ is connected as $\toba_{\bq}$ is the twisted tensor product of 
the Nichols algebras of the connected components of $\bq$.
Under suitable hypotheses, a matrix $\bq$ gives rise to a generalized root system \cite{HY-groupoid}; 
if $\dim \toba_{\bq} < \infty$, then  $\bq$ has a finite root system. 
The classification of all $\bq$ 
with finite root system and connected Dynkin diagram was given in \cite{H-classif}; this contains
the classification of the finite-dimensional Nichols algebras of diagonal type.
The Nichols algebras $\toba_{\bq}$ with $\bq$ in the list of \cite{H-classif} have finite $\GK$. It was conjectured
that these are all.

\begin{conjecture}\label{conj:AAH} \cite [Conjecture 1.5] {AAH-triang} The root system 
of a Nichols algebra of  diagonal type with finite GK-dimension is finite. 
\end{conjecture}

The conjecture holds when $\theta \leq 3$,  when $\bq$ is of  Cartan type, or when $\bq$ is generic, 
see \cite{AAH-diag},  \cite{angiono-garcia2} and references therein.

The defining relations of the Nichols algebras $\toba_{\bq}$ 
with $\bq$ in the list of \cite{H-classif} appear in \cite{angiono-convex,angiono-diagonal}. See the survey \cite{AA-diag}.

\medbreak
We shall apply several times the following result. The first three items are well-known, the last follows from 
\cite[Theorems 1, 2, 4.1]{AAH-diag}. The case $\theta =1$ is covered by Example \ref{exa:YD-pair}.

\begin{lemma}\label{lema:diagonal-q-constante} Assume that $\theta \geq 2$.
Given $q \in \kut$, 
let $V$ be a braided vector space of diagonal type with matrix $\bq = (q_{ij})_{i, j \in \I}$ where $q_{ij} = q$ for all $i,j\in \I$; thus locally the Dynkin diagram is 
$\xymatrix{\underset{ i }{\overset{q}{\circ}} \ar  @{-}[r]^{q^2}   & \underset{ j }{\overset{ q }{\circ} }}$ 
for all $i \neq j$. Then
\begin{enumerate}[leftmargin=*,label=\rm{(\roman*)}]
\item\label{item:diagonal-q-constante1} If $q = 1$, then $\toba(V) \simeq S(V)$.

\item\label{item:diagonal-q-constante2} If $q = -1$, then $\toba(V) \simeq \Lambda(V)$.

\item\label{item:diagonal-q-constante3}  If $q \in \G'_3$ and $\dim V =2$, then $\toba(V)$ is of Cartan type $A_2$ and has dimension 27. 

\item\label{item:diagonal-q-constante4} Otherwise $\GK \toba(V) = \infty$. \qed
\end{enumerate}
\end{lemma}

\medbreak
The next result will  be useful too. For the first two items we assume Conjecture \ref{conj:AAH}
and go through the list of \cite{H-classif}. Here by cycle we mean a closed path.
For \ref{item:rossito}, see \cite[Lemma 2.8] {AAH-triang} inspired by \cite{R quantum groups}.

\begin{lemma}\label{lema:diagonal-ciclos}
Let $V$ be a braided vector space of diagonal type such that its Dynkin diagram contains 

\begin{enumerate}[leftmargin=*,label=\rm{(\roman*)}]
\item\label{item:cycle>3} either an $N$-cycle, with $N > 3$;

\item\label{item:cycle=3} or else a $3$-cycle, with no vertex labelled
by $-1$,

\item\label{item:rossito} or else a sub-diagram of the form  
$\xymatrix{ \underset{ \ }{\overset{q}{\circ}} \ar  @{-}[r]^{r}   & \underset{ \ }{\overset{ 1 }{\circ} }}$, $r\neq 1$.
\end{enumerate}
Then $\GK \toba(V) = \infty$. \qed
\end{lemma}

\begin{remark}\label{remark:3-cycle} The only $3$-cycles with all vertices labelled by $-1$ in \cite{H-classif} are 
\begin{align*}
&\raisebox{0.5cm}{\xymatrix@C-4pt@R-10pt{
&\overset{-1}{\circ} \ar  @{-}[dl]_{s} \ar  @{-}[dr]^{r} & \\
\overset{-1}{\circ} \ar  @{-}[rr]^{q}& & \overset{-1}{\circ}}}&
&\text{where $q,r,s\neq 1$, $qrs=1$;}
\end{align*} 
These are of type $D(2, 1; \alpha)$,  \cite[\S 5.4]{AA-diag}.
When $q=r=s  \in \G'_{3}$, $\dim \toba_{\bq} = 2^4 3^3$.
Furthermore this diagram can not be embedded in a Dynkin diagram of rank 4 with finite root system, 
see \cite[Table 4]{H-classif}.
\end{remark}

\begin{remark}\label{remark:2-cycle} The only Dynkin diagrams of rank 2 in the list of \cite{H-classif}
 with both vertices labelled by the same $q \neq 1$ are
 
\medbreak
\noindent
$\xymatrix@C-10pt{
\overset{q}{\circ} \ar  @{-}[rr]^{q^{-1}}& & \overset{q}{\circ}}$, of type $A_2$;

\medbreak
\noindent$\xymatrix@C-10pt{
\overset{q}{\circ} \ar  @{-}[rr]^{\zeta}& & \overset{q}{\circ}}$, 
where $\zeta \in \G'_{12}$ and $q = -\zeta^{2}$; 
of type $\mathfrak{ufo}(8)$, see \cite[\S 10.8]{AA-diag}.
\end{remark}

\begin{definition}\label{remark:ppal-realiz-odd-order}
Let $(V,c^{\bq})$ be a braided vector space of diagonal type. A principal realization of $(V,c^{\bq})$ over 
a Hopf algebra $\cH$ is a family $(g_i, \chi_i )_{ i \in \I}$ of YD-pairs 
such that $q_{ij} = \chi_j(g_i)$ for all $i, j$. In this case $V = \bigoplus_i \ku^{\chi_i}_{g_i} \in \yd{\cH}$.
\end{definition}

\medbreak
Let $G$  be a finite group  of odd order. By inspection of the list in \cite{H-classif}, we see that
a matrix $\bq$ with finite  root system and connected Dynkin diagram
could have a  principal realization  over the the group algebra $\ku G$
only when either it is of Cartan type, or else its Dynkin diagram is one of
\begin{align}\label{eq:dynkin-br(2,a)}
&\xymatrix{  \overset{\omega}{\underset{\ }{\circ}} \ar  @{-}[r]^{q^{-1}}  &
\overset{q}{\underset{\ }{\circ}}},
& &\xymatrix{  \overset{\omega}{\underset{\ }{\circ}} \ar  @{-}[rr]^{\omega^2q}  &&
\overset{\omega q^{-1}}{\underset{\ }{\circ}}},
\\
\label{eq:dynkin-br(3)}
&\xymatrix{ \overset{\zeta}{\underset{\ }{\circ}}\ar  @{-}[r]^{\ztu}  &
\overset{\zeta}{\underset{\ }{\circ}} \ar  @{-}[r]^{\ztu}  & \overset{\zeta^{-3}}{\underset{\ }{\circ}}}
& &\xymatrix{  \overset{\zeta}{\underset{\ }{\circ}}\ar  @{-}[r]^{\ztu}  &
\overset{\zeta^{-4}}{\underset{\ }{\circ}} \ar  @{-}[r]^{\zeta^{4}}  & \overset{\zeta^{-3}}{\underset{\ }{\circ}}},
\end{align}
where the order of $q$ divides $\vert G\vert$ and is $>3$, $\omega\in \G_{3}'$ and $\zeta\in \G_{9}'$.
In particular, if $3$ does not divide $\vert G\vert$, then $\bq$ could admit a principal realization only 
if it is of Cartan type. See \cite[\S 7.2, \S 7.3]{AA-diag} for information on \eqref{eq:dynkin-br(2,a)} of type $\mathbf{br}(2)$, 
resp. \eqref{eq:dynkin-br(3)} of type $\mathbf{br}(3)$.

\subsection{Yetter-Drinfeld modules over groups}\label{subsection:ydg}
Let $G$ be a group. 
Recall that $\mathcal O_x$ denotes the conjugacy class and $G^x$ the centralizer of $x\in G$. 
For any $y \in \Oc_x$ we
fix $g_y \in G$ such that $g_y \trid x   = y$. 
Then for $h\in G$ and $y \in \Oc_x$ 
\begin{align}\label{eq:thy}
t_{h,y} &:= g_{h\trid y}^{-1}hg_y \in G^x . 
\end{align}
A Yetter-Drinfeld module $M \in \yd{\ku G}$ is just a $G$-graded vector space $M = \oplus_{g \in G} M_{g}$
provided with a linear action of $G$ such that
\begin{align*} 
h \cdot M_{g} &= M_{h\trid g}, & h,g &\in G. 
\end{align*}
In such case, the support of $M$ is $\supp M  = \{g\in G: M_g \neq 0\}$, which is
a disjoint union of conjugacy classes. If $v \in M_{g}$, then $\deg v := g$.

\medbreak We next  describe $\Irr \yd{\ku G}$ and $\Indec \yd{\ku G}$. First we consider $x \in G$ and 
a  representation $\rho: G^x\to GL(W)$. We set 
\begin{align*}
M(x, W) := \Ind_{G^x}^{G} W
\simeq\ku \mathcal O \otimes W \simeq \bigoplus_{y \in \Oc_x} g_y
\otimes W.
\end{align*}
As is known,  $M(x, W)$ belongs to $\yd{\ku G}$ with  action and  grading 
\begin{align*}
h\cdot(g_y\otimes w) & =g_{h\trid y}\otimes t_{h,y} \cdot w,&
\deg (g_y \otimes w) &= y,
\end{align*}
where $t_{h,y}$ is given by \eqref{eq:thy}.
We also use the notation $M(x, \rho) = M(x, W)$, and accordingly 
\begin{align*}
\toba(x, W) := \toba\left(M(x, \rho)\right) =: \toba(x, \rho).
\end{align*}
For brevity, we set $g_y w = g_y\otimes w$.
Then the braiding of $M(x, W)$ is given by
\begin{align}\label{eq:braiding-ydG}
c(g_z u\otimes g_y w) &= g_{z\trid y} (t_{z,y}\cdot w)\otimes g_z u,& z,y \in \Oc_x,& u,w \in W.
\end{align}

\begin{example}\label{exa:YD-pair-group}
A YD-pair over $\ku G$ is just a pair
$(g, \chi) \in Z(G) \times \widehat G$; then $\ku_g^{\chi} \simeq M(\Oc_g, \chi)$, 
see Example \ref{exa:YD-pair}.
If e.g. $g \in Z(G) \cap [G,G]$ and $\chi \in \widehat G$,
then $\toba(g, \chi) \simeq S(W)$ where $\dim W = 1$, thus $\GK \toba(g, \chi) = 1$. \qed
\end{example}

\begin{prop}\label{prop:simpleYD}
$\Irr \yd{\ku G}$ is parametrized by pairs $(\Oc, W)$ where $\Oc$ is a conjugacy class and $W \in \Irr G^x$
for a fixed choice of $x \in \Oc$.
Similarly $\Indec \yd{\ku G}$ is parametrized by pairs $(\Oc, W)$ where now $W \in \Indec G^x$.
\end{prop}
We sketch a proof of this well-known result for completeness, as we have not found a reference for the case when $G$ is infinite.

\pf
If $M$ is indecomposable, then necessarily $\supp M  = \Oc_x$ for some $x \in G$
that we fix.
Thus $M = \oplus_{y \in \Oc_x} M_y$ and $G^x$ acts on $M_x$. Let $N_x$ be a $G^x$-submodule of $M_x$.
For each $y \in \Oc_x$ choose $g_y \in G$ such that $g_y \trid x = y$.
Then $g_y \cdot N_x \subset M_y$ and $N := \oplus_{y \in \Oc_x} g_y \cdot N_x$
is a Yetter-Drinfeld submodule of $M$. Thus, if $M$ is indecomposable, respectively simple, then so is 
$M_x$ as $G^x$-module, and $M \simeq M(x, M_x)$.
The converse is proved similarly.
\epf

\subsection{Racks}\label{subsection:typeC} Nichols algebras over groups are studied systematically through racks. 
We refer to \cite{AG-adv} for an exposition on racks and Nichols algebras over groups  and \cite{A-leyva} for more recent results. Here we  collect some material needed in this paper.
A rack is a non-empty set $X$ with a self-distributive
operation $\trid: X \times X \to X$ such that $\varphi_x := x\trid \underline{\ }$ is 
bijective for every $x \in X$. The main examples are subsets of groups stable under conjugation.
All racks here are assumed to be  subracks of groups. 
A rack $X$ is \emph{abelian} if $x \trid y = y$ for all $x,y\in X$.
The inner group of a rack $X$ is the subgroup $\Inn X$ of the group $\Aut X$ of rack automorphisms generated by
$\varphi_x$ for all $x \in X$. 

\begin{lemma}\label{lema:Inn-pi} \cite[Lemma 1.8]{AG-adv}
A surjective morphism of racks $\pi: X \to Y$ extends to a surjective  morphism of groups
$\Inn \pi: \Inn  X \to \Inn  Y$.
\end{lemma}

\pf Given $x \in X$, define $\Inn \pi(\varphi_x) = \varphi_{\pi(x)}$. If  $z \in X$ satisfies 
$\varphi_x = \varphi_z$, then $\varphi_{\pi(x)} (\pi(y)) = \pi(x) \trid \pi(y) =\pi(x \trid y)=\pi(z \trid y) 
= \varphi_{\pi(z)} (\pi(y))$. Since $\pi$ is surjective, then $\varphi_{\pi(x)} = \varphi_{\pi(z)}$, i. e. 
$\Inn \pi$ is well-defined. Consider next
\begin{align*}
G = \{\sigma \in \Sb_X: \exists \nu \in \Sb_Y \text{ such that } \pi \sigma = \nu \pi \}.
\end{align*}
Clearly such $\nu$ is unique, $G \leq \Sb_X$ and $\sigma \mapsto \nu$ is a  morphism of groups.
Thus $\Inn \pi$ is a   morphism of groups which is surjective because $\pi$ is so.
\epf

The following statement is a consequence of  \cite[Theorem 2.1]{HV}.
\begin{theorem}\label{th:HV} 
Let $G$ be a finite non-abelian group and $V$ and $W$ be two  simple Yetter-Drinfeld
modules over $G$ such that $G$ is generated by the support of $V\oplus W$, $\dim V \le \dim W$ and 
\begin{align}
{c^2}_{\vert V\otimes W} \neq \id_{V\otimes W}.
\end{align}
If $\dim \toba (V\oplus W) < \infty$, then
$(\dim V, \dim W)$ belongs to
\begin{align}\label{eq:rank2-set}
\{(1,3), (1,4), (2,2), (2,3), (2,4)\}. 
\end{align}
When $\dim W =3$ in \eqref{eq:rank2-set}, $W$ is the braided vector space associated to the transpositions in 
$\mathbb{S}_3$ with the cocycle $-1$, which is not of diagonal type.
\qed\end{theorem}

Recall from \cite[Definition 2.3]{ACG-III} that a finite rack
$X$ is \emph{of type C}  when there are  a decomposable subrack
$Y = R\coprod S$  and elements $r\in R$, $s\in S$ such that
\begin{align}
\label{eq:typeC-rack-inequality}
&r \triangleright s \neq s \qquad (\text{hence } s \trid r \neq r),
\\ \label{eq:typeC-rack-indecomposable}
&R = \Oc^{\Inn Y}_{r}, \qquad S = \Oc^{\Inn Y}_{s},
\\
\label{eq:typeC-rack-dimension} &\min \{\vert R \vert, \vert S \vert \}  > 2 \quad \text{ or } \quad
\max \{\vert R \vert, \vert S \vert \}  > 4.
\end{align}

\begin{theorem}\label{th:typeC} \cite[Theorem 2.9]{ACG-III}
A finite rack of type C collapses, that is
$\dim\toba(\Oc, \bq) = \infty$ for every finite faithful 2-cocycle $\bq$. \qed
\end{theorem}

The proof of Theorem \ref{th:typeC} relies on \cite[Theorem 2.1]{HV}, 
which in turn depends on the notion of Weyl groupoid \cite{HY-groupoid}.
For some of the arguments below we need the validity of the following conjecture; the adaptation 
of the proof of \cite{HV} does not appear to be straightforward.

\begin{conjecture}\label{conjecture:typeC-X-finite-GK}  Let $X$ be a finite rack  of type C.
Then $\GK \toba(\Oc, \bq) = \infty$ for every faithful 2-cocycle $\bq$.
\end{conjecture}

\begin{example}\label{exa:yd-product-groups}\label{remark:product-typeC}
 Assume that $G = G_1 \times G_2$.
If $x = (x_1, x_2) \in G$, then
\begin{align} \label{eq:cclas-product}
\Oc_x &= \Oc_{x_1} \times \Oc_{x_2}, &G^x &= G^{x_1} \times G^{x_2}, & \Irr G^x &\simeq \Irr G^{x_1} \times \Irr G^{x_2}.
\end{align}
Thus if $W \simeq W_1 \otimes W_2$ is a simple $G^x$-module, then
\begin{align}\label{eq:yd-product-groups}
M(x, W) & \simeq M(x_1, W_1) \otimes M(x_2, W_2).
\end{align}
Correspondingly, the tensor product of two braided vector spaces $(V, c_V)$ $(W, c_W)$
is $(V \otimes W, c_{V\otimes W})$ where $c_{V\otimes W}: V\otimes W \otimes V\otimes W \to
V\otimes W \otimes V\otimes W$ is defined by 
\begin{align*}
c_{V\otimes W} := \left(\id \otimes \tau_{V, W} \otimes \id\right)\left(c_{V} \otimes c_{W}\right) \left(\id \otimes \tau_{W, V} \otimes \id\right).
\end{align*}
Observe that there is no clear relation between the Nichols algebras
$\toba(x, W)$, $\toba(x_1, W_1)$ and $\toba(x_2, W_2)$.
For instance if all three modules in \eqref{eq:yd-product-groups} have dimension one, then the braidings
are given respectively by $q,q_1$ and $q_2$ with $q = q_1q_2$ so that any of them could $1$ 
with the other two being non-trivial roots of 1.
But the criterium of type C propagates in this setting. 

\smallbreak
Namely, let
$X = X_1 \times X_2$ be a direct product of racks. If either $X_1$ or $X_2$ is of type C, then so is $X$. 
Indeed, let $Y_1 = R_1\coprod S_1$ be a subrack of $X_1$ and elements $r_1\in R_1$, $s_1\in S_1$ satisfying \eqref{eq:typeC-rack-inequality}, \eqref{eq:typeC-rack-indecomposable} and \eqref{eq:typeC-rack-dimension}.
Pick any $x_2 \in X_2$ and set $Y = Y_1\times \{x_2\}$, $R = R_1\times \{x_2\}$, $S = S_1\times \{x_2\}$,
$r = (r_1, x_2)$, $s = (s_1, x_2)$. Then \eqref{eq:typeC-rack-inequality}, \eqref{eq:typeC-rack-indecomposable} and \eqref{eq:typeC-rack-dimension} hold for them.
\end{example}

\section{Hopf algebras and conjugacy classes}\label{section:hopf-nilpotent}

Recall that the upper 
central series  of a group $G$ is the sequence of subgroups $e = Z_0 \lhd Z_1 \lhd \cdots \lhd Z_i \lhd \cdots$, where
\begin{align*}
Z_{n+1} &= Z_{n+1}(G) = \{x\in G:[x,G] \leq Z_n\};
\end{align*}
$G$ is nilpotent iff the upper  centralizer series stabilizes in $G$ \cite[Th. 2.2]{CMZ}.

\subsection{Conjugacy classes in finite nilpotent groups}\label{subsection:conjugacy-finite-nilpotent}
Here is our first basic result on Nichols algebras over finite nilpotent groups.

\begin{theorem}\label{thm:odd-order-abelian-or-typeC}
Let $\Oc$ be a conjugacy class in a finite nilpotent group $G$ of odd order. Then $\Oc$
is either of type C or else an abelian rack.
\end{theorem}

There are examples of conjugacy classes of finite nilpotent groups that are of  type C, see page \pageref{section:nichols-triangular3}. 

\pf It is well-known that a finite group is nilpotent if and only if it is isomorphic to the product of its Sylow subgroups,
see e.~g. \cite[Theorem 2.13]{CMZ}. Hence, by  Example \ref{remark:product-typeC}, 
we may assume that $G$ is a $p$-group with $p$ an odd prime.
Let us assume that $\Oc$ is  not abelian. That is, there exist $r,s \in \Oc$ such that $r\trid s \neq s$ 
(and then $s\trid r \neq r$). Let $H = \langle r,s \rangle \leq G$.

If $R:=\Oc_r^{H}\neq S:=\Oc_s^{H}$, 
then  $Y:=R\coprod S$ is a decomposable subrack of $\Oc$
that satisfies \eqref{eq:typeC-rack-indecomposable}, because  $H = \langle Y\rangle$ so 
$R=\Oc_r^{H}=\Oc_r^{\Inn Y}$ and $\Oc_s^{H}=\Oc_s^{\Inn Y}$.  
Now \eqref{eq:typeC-rack-inequality} holds by assumption. 
Evidently $p$ divides both $\vert R \vert$ and $\vert S \vert$. 
Since  $p \geq 3$, \eqref{eq:typeC-rack-dimension} holds and $\Oc$ is of type C.

\medbreak
Next suppose that $s \in \Oc_r^{H} = Y$. Then $H = \langle Y \rangle$, hence $Y = \Oc_r^{H} = \Oc_r^{\Inn Y}$ 
is indecomposable by \cite[Lemma 1.15]{AG-adv}. Also  $\Inn Y \simeq H/Z(H)$ by \cite[Lemma 1.9]{AG-adv}, hence 
$\Inn Y$ is a $p$-group. Now by a routine recursive argument, there exists a surjective morphism of racks
$\pi: Y \to Z$ where $Z$ is simple. Considering the surjective morphism of groups  $\Inn \pi: \Inn Y \to \Inn Z$ 
given by Lemma \ref{lema:Inn-pi}, we see that  $\Inn Z$ is a $p$-group, 
so in particular $\vert Z\vert$ should be a power of $p$.
But then $\Inn Z$ could not be a $p$-group being a semidirect product of a $p$-group with a group of order not divisible by $p$, 
see \cite[Proposition 3.2 and Theorem 3.7]{AG-adv}, and also
 the discussion at the end of page 204 and the beginning of  page 205 in \cite{AG-adv}. 
\epf

Note that this is not a statement on racks with $p^n$ elements, $p$ an odd prime but on conjugacy classes of $p$-groups.

\medbreak
When $G$ is a $2$-group, Theorem \ref{thm:odd-order-abelian-or-typeC} is no longer true.
Indeed, the group $\D_4 = \langle x,y \vert x^2 = e =y^4, xyx = y^3\rangle$ admits a a finite-dimensional Nichols algebra
$\toba(V)$ where $\supp V = Y = \Oc_x \coprod \Oc_{xy}$
which is neither abelian nor of type C. The rack $Y$ can be realized as conjugacy class in the group $D_4 \rtimes\Z/4$
determined by the automorphism $\sigma: \D_4 \to \D_4$,  $\sigma(x) = xy$, $\sigma(y) = y$.

\subsection{Hopf algebras and infinite conjugacy classes}\label{section:hopf}
Let $\cH$ be a Hopf algebra.
If $R$ is a Hopf algebra in $\yd{\cH}$ then the \emph{bosonization}  (or biproduct) $R \# \cH$ is the Hopf algebra  with
underlying vector space $R\otimes \cH$ and structure as in \cite[Section 11.6]{Rad-libro}.

\begin{remark}\label{rem:fingen-discrete}
Recall that an affine algebra is a finitely generated one.
Let $\Oc$ be a conjugacy class in a finitely generated group $G$.
Let $x \in \Oc$ and let $W \in \Rep {G^x}$ be a finitely generated module.
Set $M = M(\Oc, W)$. Then $T(M)\# \ku G$ is affine, hence so is $\toba(M)\# \ku G$. 
\end{remark}

\pf Let $(g_i)_{i\in I}$ be a family of generators of $G$ and $(w_j)_{j \in J}$ be a family of generators of $W$.
Then the $g_i$'s together with the $w_j$'s generate the algebra $T(M)\# \ku G$.
\epf

Let $G$ be a finitely generated group and let $M \in \yd{\ku G}$ such that the action of $G$ is locally finite.
By \cite[Lemma 2.3.1]{AAH-triang}, we have 
\begin{align}\label{eq:GK-bos}
\GK \toba(M) \# \ku G \le \GK \toba(M) + \GK \ku G.
\end{align} 
Furthermore, if $\dim M < \infty$, then the equality holds in \eqref{eq:GK-bos}.
Our second basic result, Theorem \ref{thm:infinite-conjugacy-class} (inspired by \cite[Example 2.3.3]{AAH-triang}), 
roughly states that $\GK \toba(M) \# \ku G = \infty$ if the support of $M$ is an infinite conjugacy class,
even if $\GK \toba(M)$ is finite, in sharp contrast with \eqref{eq:GK-bos}.
We start by a theorem of Malcev needed for our approach.

\begin{theorem}\label{theorem:malcev-finite-index} \cite[Theorems 2.23, 2.24]{CMZ} (Malcev) 
Let $G$ be a finitely generated nilpotent group and let $H \leq G$.
Assume that there exists a finite set $X$ of generators such that
for any $g \in X$ there exists a positive integer $n$ such that $g^n \in H$.
Then the index of $H$ in $G$ is finite. 

Furthermore if for any $g \in X$ the integer $n$ is a power of a fixed prime $p$, 
then $[G:H]$ is also a power of $p$. \qed
\end{theorem}

Actually, the last claim holds more generally if $n$ is a $\varpi$-number, where $\varpi$ is a fixed set of primes.

\begin{coro}\label{coro:malcev-finite-index}  
Let $G$ be a finitely generated nilpotent-by-finite group and let $H \leq G$.
Assume that for every $g \in G$ there exists a positive integer $n$ such that $g^n \in H$.
Then the index of $H$ in $G$ is finite.
\end{coro}

\pf Let $N \leq G$ be nilpotent of finite index; $N$ is finitely generated by \cite[1.6.11]{robinson}.
If $g \in N$, then there exists $n \in \N$ such that $g^n \in H \cap N$; thus $[N:H\cap N]$ is finite by Theorem \ref{theorem:malcev-finite-index} and so is $[G: H \cap N] = [G:N][N:H\cap N]$.  But $[G: H \cap N] = [G:H][H:H\cap N]$,
so $[G:H]$ is finite.
\epf

\begin{coro}\label{coro:nilpotent-by-finite-infinite-orbit}  
Let $G$ be a finitely generated nilpotent-by-finite group, let $\Oc \subset G$ be an infinite conjugacy class
and pick $x \in \Oc$.
Then there exists $g \in G$ such that $g^n \trid x \neq x$ for all $n \in \N$.
\end{coro}

\pf
If for every $g \in G$ there exists $n  \in \N$ such that $g^n \trid x = x$, then 
$\vert \Oc \vert = [G:G^x]$ is finite by Corollary \ref{coro:malcev-finite-index}.
\epf

\begin{theorem}\label{thm:infinite-conjugacy-class}
Let $G$ be a finitely generated group and $M\in \yd{\ku G}$ such that $\Oc = \supp M$ is an infinite conjugacy class.
Then 
\begin{align*}
\GK \toba(M) \# \ku G = \infty.
\end{align*}
\end{theorem}

\pf  By Gromov's Theorem we may assume that $G$ is nilpotent-by-finite. 
Let $x \in \Oc$. Since $M_x$ is the union of its finitely generated $G^x$-submodules,
we may assume that it is finitely generated.
Let $S$ be a finite set of generators of the group $G$
and let $F$ be a finite set of generators of the $G^x$-module $M_x$.
Pick $g \in G$ such that $x_n = g^n \trid x \neq x$ for all $n \in \N$, which exists by Corollary \ref{coro:nilpotent-by-finite-infinite-orbit}, and $m_0 \in F \backslash 0$.
Let $V = \langle 1, S, g^{\pm 1}, F\rangle$, a set of generators of $\toba(M) \# \ku G$. 
For $n \in \N$ we set
$m_n = g^n m_0g^{-n}\in V^{2n+1}$. 
Given $s \in \I_{0,n}$, let
\begin{align*}
\varLambda_{n,s} &= \{m_{i_1} \cdots m_{i_s}: 1 \leq i_{1} < \cdots < i_{s} \leq n \},
&
\varLambda_{n} &= \bigcup_{s \in \I_{0,n}} \varLambda_{n,s}.
\end{align*}
We claim that (i) $\vert \varLambda_{n} \vert = 2^{n}$ (exercise), (ii)
$\varLambda_{n}  \subset V^{(n+1)^{2}}$, and (iii) $\varLambda_{n}$  is a linearly independent set.
For (ii), just observe that 
\begin{align*}
\varLambda_{1} &= \{1, m_1\} \subset  V^{3}\subset  V^{4},&
\varLambda_{n} & \subset  \varLambda_{n-1} \{1, m_n\} \subset  V^{n^{2}}V^{2n+1} \subset V^{(n+1)^{2}}.
\end{align*}
For (iii), it is enough to prove that $\varLambda_{n,s} \subset \toba^s(M)$ is a linearly independent set.
Clearly, $m_n \in M_{x_n}$ for all $n \in \N_0$. The claim (iii) is a particular case of 

\setcounter{claim}{3}
\renewcommand{\theclaim}{(\roman{claim})}
\begin{claim}\label{claim4} Given $\widetilde{m}_i \in M_{x_i} \backslash 0$, $i \in \I_{n}$, the subset
\begin{align*}
\widetilde{\varLambda}_{n,s} = \{\widetilde{m}_{i_1} \cdots \widetilde{m}_{i_s}: 1 \leq i_{1} < \cdots < i_{s} \leq n \}
\end{align*}
of $\toba(M)$ is linearly independent.
\end{claim}

\noindent \emph{Proof of Claim \emph{\ref{claim4}}.} By induction on $s$ and $n$.
The elements $x_i$ are all different by our choice of $g$, thus the case $s=1$ follows. 
After completing appropriately the family $(\widetilde{m}_i)$ to a homogeneous basis of $M$,
we know that there exist skew-derivations $\partial_j: \toba(M) \to \toba(M)$, $j \in \I_{n}$,
such that
\begin{align*}
\partial_j (uv) &= \partial_j (u) (x_j\cdot v) + u\, \partial_j (v),&
\partial_j (\widetilde{m}_{i}) &= \delta_{i,j}.
\end{align*}
Set $\bI_s = \{\bi = (i_1,\cdots, i_s): 1 \leq i_{1} < \cdots < i_{s} \leq n \}$, 
$\widetilde{m}_{\bi} = \widetilde{m}_{i_1} \cdots \widetilde{m}_{i_s}$. Then
\begin{align*}
\partial_j (\widetilde{m}_{\bi}) &=  \begin{cases}
0, & \text{if } j \notin \{i_{1}, \cdots, i_{s}\},
\\ \widetilde{m}_{i_1} \cdots \widetilde{m}_{i_{h-1}} (x_j \cdot \widetilde{m}_{i_{h+1}}) \cdots (x_j \cdot \widetilde{m}_{i_{s}})
& \text{if } j = i_{h}, \ h \in \I_s.
\end{cases}
\end{align*}
Thus we consider the map $\psi_j$ from $\bI_{s;j} :=\{\bi =  (i_1,\cdots, i_s) \in \bI_s: \exists h, i_h = j \}$
to $\bI_{s-1; \neg j} :=\{\bi =  (i_1,\cdots, i_{s-1}) \in \bI_{s-1}: \nexists k, i_k = j \}$ given by
\begin{align*}
\bi &\mapsto (i_1,\cdots, i_{h-1}, i_{h+1}, \cdots i_s).
\end{align*}
It is easy to see that this map $\psi_j$ is bijective. 
Let now $\lambda_{\bi}$, $\bi\in \bI_s$, be a family of scalars such that
$\sum_{\bi\in \bI_s} \lambda_{\bi} \widetilde{m}_{\bi} = 0$. Then for any $j \in \I_n$
\begin{align*}
0 & = \partial_j \big(\sum_{\bi\in \bI_s} \lambda_{\bi} \widetilde{m}_{\bi} \big)
= \sum_{\bi\in \bI_{s;j}} \lambda_{\bi} \widetilde{m}_{i_1} \cdots \widetilde{m}_{i_{h-1}} (x_j \cdot \widetilde{m}_{i_{h+1}}) \cdots (x_j \cdot \widetilde{m}_{i_{s}})
\end{align*}
If $j=n$, then
\begin{align*}
0 &  = \sum_{\bi\in \bI_{s;n}} \lambda_{\bi} \widetilde{m}_{i_1} \cdots \widetilde{m}_{i_{s-1}} \overset{\star}{\implies}
\lambda_{\bi} = 0 \ \forall   \bi\in \bI_{s;n} \implies
\\
0 &  = \sum_{\bi\in \bI_{s; \neg n}} \lambda_{\bi} \widetilde{m}_{\bi} = \sum_{\bi\in \bI_{s; n-1}} \lambda_{\bi} \widetilde{m}_{\bi} \overset{\diamond}{\implies} \lambda_{\bi} = 0 \ \forall   \bi\in \bI_{s;\neg n},
\end{align*}
where $\star$ is by the inductive hypothesis on $s$ and $\diamond$ is by the inductive hypothesis on $n$. 
Claim {\ref{claim4}} is proved.

By (i), (ii) and (iii), $2^{n} \leq \dim V^{(n+1)^{2}}$ for all $n$; the Theorem follows.
\epf

\subsection{Finitely generated torsion-free nilpotent groups}\label{subsection:nichols-nilpotent-general-main}
Recall that the FC-centre $FC(G)$ of a group $G$
is the union of all finite conjugacy classes of $G$;
$FC(G)$ is a characteristic subgroup of $G$  containing $Z(G)$ \cite{baer}.
The following result is folklore; see for instance \cite{eckhardt-gillaspy} for a different proof.

\begin{lemma}\label{lema:noncentral-conjugacyclasses}
If $G$ is a finitely generated torsion-free nilpotent group, then  every non-central conjugacy class is infinite.
\end{lemma}

\pf Since $G$ is finitely generated nilpotent, so is $FC(G)$. By \cite[Theorem 1.6]{tomkinson},
 $[FC(G), FC(G)]$ is finite. Since $G$ is torsion-free,  
$[FC(G), FC(G)] = e$ i.e. $FC(G) = Z(FC(G)) = Z(G)$, the last equality by \cite[Lemma 2.2]{eckhardt-gillaspy}.
\epf

\subsection{Finite conjugacy classes in nilpotent groups}\label{subsection:odd-torsion}
We generalize Theorem \ref{thm:main-nichols-nilpotent-odd} to nilpotent groups whose torsion has odd order,
using a well-known result of Gruenberg.

\begin{theorem}\label{thm:gruenberg2} \cite{gruenberg} 
Let $G$ be a finitely generated  nilpotent group with torsion $T \neq e$ and $e \neq g \in G$.
Then there exists a prime $p$ that divides $\vert T \vert$,
a finite $p$-group $P$ and a morphism $\pi: G \to P$ such that $\pi(g) \neq e$.
\end{theorem}

\begin{prop}\label{prop:noncentral-conjugacyclasses-torsion-typeC}
Let $G$ be a finitely generated  nilpotent group 
whose torsion subgroup $T \neq e$ is non-trivial and has odd order. Then
a finite conjugacy class $\Oc$ of $G$ is either abelian or else of type C. 
\end{prop}

\pf Suppose that $\Oc$ is not abelian. Pick $r,s \in \Oc$ such that $r\trid s \neq s$, i.~e.
$[r,s] \neq e$. 
By Theorem \ref{thm:gruenberg2} there exist an odd prime $p$ and
an epimorphism $\pi$ from $G$ to a finite $p$-group $P$ such that 
$[\pi(r), \pi(s)] = \pi([r,s]) \neq e$.
Let $H = \langle \pi(r),\pi(s) \rangle \leq P$, $\widetilde{R} :=\Oc_{\pi(r)}^{P}$, $\widetilde{S}:=\Oc_\pi(r)^{P}$. 
Arguing as in the proof of Theorem \ref{thm:main-nichols-nilpotent-odd}, we see that
$\widetilde{Y}:= \widetilde{R}\coprod \widetilde{S}$ is a decomposable subrack of $P$. 
Let $R :=\pi^{-1}(\widetilde{R}) \cap \Oc$, $S :=\pi^{-1}(\widetilde{R}) \cap \Oc$.
Then  $Y:= R \coprod S$ is a decomposable subrack of $\Oc$. Now the elements $\pi(s)$, $\pi(r)\trid \pi(s)$,
$\pi(r)^2\trid \pi(s)$ of $\widetilde{S}$, respectively $\pi(r)$, $\pi(s)\trid \pi(r)$,
$\pi(s)^2\trid \pi(r)$ of $\widetilde{R}$, are different. Hence 
$s, r\trid s, r^2 \trid s \in S$, respectively $r, s\trid r, s^2 \trid r \in R$, are different, 
\eqref{eq:typeC-rack-dimension} holds and $\Oc$ is of type C.
\epf

Theorem \ref{thm:infinite-conjugacy-class} and Proposition \ref{prop:noncentral-conjugacyclasses-torsion-typeC}
show that the classification of pointed Hopf algebras with finite $\GK$ over 
a finitely generated  nilpotent group whose torsion has odd order goes through Nichols algebras over abelian groups.
For instance, a Hopf algebra $H$ like this is co-Frobenius if and only if  
$\gr H \simeq \toba(V)$, where $V$ is of diagonal type and $\dim \toba(V) < \infty$,
as follows from the preceding results and \cite[Theorem 1.4.2]{AAH-triang}; see \emph{loc. cit.} for details.
This last claim does not assume Conjecture \ref{conjecture:typeC-X-finite-GK}.

\bigbreak
We finish this Section with a result needed later, see Lemma \ref{lema:orbg-divides-torsion}.

\begin{lemma}\label{lema:order-conjugacy-class}
Let $G$ be a finitely generated nilpotent group with torsion $T$ and let
$\Oc = \Oc_x$ be a finite conjugacy class in $G$. 
Then $\vert \Oc \vert$ divides $\vert T \vert$.
\end{lemma}

\pf
As $G/T$ is torsion-free, $FC(G/T)=Z(G/T)$ by Lemma \ref{lema:noncentral-conjugacyclasses}. 
Hence the image of $x$ is central in $G/T$; i.~e. $[G,x]\leq T$. Let $\phi:G \to T$ be  given by
$g\mapsto [g,x]$, $g \in G$ and let $S = \phi^{-1} (T \cap Z(G))$. If $g \in G$ and $\ h \in S$, then 
\begin{align}\label{eq:phi-gh}
\phi(gh) = [gh,x] &= ghx h^{-1}g^{-1}x^{-1} = g[h,x] x g^{-1}x^{-1} = \phi(g) \phi(h).
\end{align}
Thus the restriction $\phi:S \to T \cap Z(G)$ is a homomorphism; clearly, $G^x \leq S$ and $S/G^x$ embeds into $T \cap Z(G)$. 
Hence $\vert \Oc^S_x \vert$ divides $\vert T \cap Z(G) \vert$.
Let now 
\begin{align*}
k = \min \{k\in \N_0: T\leq Z_k(G)\}. 
\end{align*}

We argue by induction on $k$ and $|\Oc|$. If $k =0$, then $G$ is torsion-free and  Lemma \ref{lema:noncentral-conjugacyclasses} applies.
If $k = 1$, then $T \leq Z(G)$,  $G = S$ and the claim follows.

\medbreak
Assume that  $k > 1$. Let  $H = G/(T\cap Z(G))$ and let $\pi: G \to H$ be the natural projection.
Given $y \in \Oc^H_{\pi(x)}$, we fix $z_y \in \Oc^G_x = \Oc$ such that $\pi(z_y) = y$ and 
$g_y \in G$ such that $z_y = g_y \trid x$. For $t \in \Oc^S_x$, fix $g_t \in S$ such that
$t = g_t \trid x$. We claim that the map 
\begin{align*}
\sigma:  \Oc^H_{\pi(x)} \times \Oc^S_x &\to \Oc, &
\sigma(y, t) &= z_yt x^{-1} = z_t [g_t, x], & y \in \Oc^H_{\pi(x)},  t \in \Oc^S_x,
\end{align*}
is a well-defined bijection. First, since $\phi(g_t) \in Z(G)$,
\begin{align*}
z_yt x^{-1} & = [g_y, x] x [g_t, x]   =  \phi(g_y) \phi(g_t) x 
\overset{\eqref{eq:phi-gh}}{=}  \phi(g_yg_t)x  =  g_yg_t \trid x \in \Oc;
\\
\sigma(y, t) &=\sigma(w, s) \implies \pi\left(z_y \phi(g_t)\right) = \pi\left(z_w \phi(g_s)\right)
\implies y =w \implies t=s.
\end{align*}
Finally, let $z \in \Oc$ and $y=\pi(z)$; then $z = z_{y}u $ where $u \in T\cap Z(G)$.
Now
\begin{align*}
z  &= z_{y} u = (g_y \trid x) (g_y \trid u)  = g_y \trid (xu) \implies xu \in \Oc;
\end{align*}
pick $g \in G$ such that $ux = xu = g \trid x$ and set $t= g \trid x$; since 
$u = (g \trid x)x^{-1} = \phi(g)$, we conclude that $g\in S$,  $t  \in \Oc^S_x$ and $\sigma(y,t) = z$.
The claim is proved. 

\medbreak
The torsion of $H$ is $T_1 = T/T\cap Z(G)$. Then
$Z_j(H) \cap T_1  = \pi \left(Z_{j + 1}(G) \cap T\right)$, $j \in \N_0$, and so 
$Z_{k-1}(H) \cap T_1  = T_1$. By the inductive hypothesis on $k-1$, $\vert\Oc^H_{\pi(x)}\vert$
divides $\vert T_1\vert$; thus $\vert\Oc\vert = \vert\Oc^H_{\pi(x)}\vert \vert \Oc^S_x\vert$ divides 
$ \vert T_1 \vert \vert T\cap Z(G) \vert =  \vert T \vert$.  
\epf

\section{Nichols algebras}\label{section:nichols-nilpotent}

\subsection{First remarks: $\dim W \geq 2$}\label{subsection:first-remarks}
In this subsection, we fix 
\begin{itemize}[leftmargin=*,label=$\circ$]
\item a group $G$, a finite conjugacy class $\Oc$ in $G$, $x \in \Oc$.

\medbreak
\item $W \in \Irr G^x$ with representation $\rho: G^x \to GL(W)$.  We assume that $\dim W$ is countable;
this is the case if $G$ is finitely generated. 
By the Schur Lemma, aka Dixmier's Lemma \cite[0.5.2]{wallach},
there exists  $\eta \in \widehat{Z(G^x)}$ 
implementing the action of $Z(G^x)$ on $W$.

\medbreak
\item   $\chi \in \widehat{G^x}$; set $q :=\chi\left(x\right)$. 
\end{itemize}

\begin{remark}\label{rem:change-of-point}
Pick $y \in \Oc$ and $g_y\in G$  such that $g_y \trid x =y$. Then $G^y = g_y \trid G^x$ and
$\rho^y: G^y \to GL(W)$, $\rho^y (g) = \rho(g_y^{-1} \trid g)$ is an irreducible representation
of $G^y$. Hence $M(\Oc, \rho) \simeq M(\Oc, \rho^y)$ and $\eta(x) = \eta^y(y)$. 
Similarly, $M(\Oc, \chi) \simeq M(\Oc, \chi^y)$ and $q =\chi^y\left(y\right)$.
\end{remark}

\medbreak
We start with an argument going back to \cite[3.1]{Grana}, based on Lemma \ref{lema:diagonal-q-constante}.

\begin{lemma}\label{lema:old-matiasg-GK}
Let $W \in \Irr G^x$ be as above. Assume that $\dim W \geq 2$. 
Let $\Mtt = \ku x \otimes W$, a braided subspace of $M(\Oc, W)$
with $\dim \Mtt = \dim W$. 
Then $\toba(\Mtt)$ is a braided Hopf subalgebra of $\toba(\Oc, W)$ and we have

\begin{enumerate}[leftmargin=*,label=\rm{(\roman*)}]
\item\label{item:old-matiasg-GK1} If $\eta(x) = 1$, then $\toba(\Mtt) \simeq S(\Mtt)$.
Thus $\GK \toba(\Oc, W) = \infty$ implies that $\dim W < \infty$.

\medbreak
\item\label{item:old-matiasg-GK2} If $\eta(x) = -1$, then $\toba(\Mtt) \simeq \Lambda(\Mtt)$.

\medbreak
\item\label{item:old-matiasg-GK3}  If $\eta(x) \in \G'_3$ and $\dim W =2$, then $\toba(\Mtt)$ is of Cartan type $A_2$ and has dimension 27. 

\medbreak
\item\label{item:old-matiasg-GK4} In any other case, $\GK \toba(\Oc, W) = \infty$. 
\end{enumerate}
\end{lemma}

\pf Choose $g_x = x$, thus $t_{x,x} = x$, cf. \eqref{eq:thy}. 
Fix a basis $(w_i)_{i\in I}$ of $W$,
so that the symbols $xw_i$ form a basis of $\Mtt$ and its braiding  is given by
\begin{align*}
c(xw_i\otimes xw_j) &= \eta(x) \, xw_j\otimes xw_i,& i,j \in I.
\end{align*}
Then $\GK \toba(\Mtt)$ can be read off from Lemma \ref{lema:diagonal-q-constante}.
\epf

We next generalize \cite[3.5]{YZ}.

\begin{lemma}\label{lema:old-matiasg} 
Let $G$ be a  finite group of odd order and let $W \in \Irr G^x$ as above. 

\begin{enumerate}[leftmargin=*,label=\rm{(\alph*)}]
\item\label{item:old-matiasg1} If $\dim W \geq 2$ and
$\GK \toba(\Oc, W)$ is finite, then 
$\eta(x) = 1$ and consequently $\GK \toba(\Oc, W) >0$.

\medbreak
\item\label{item:old-matiasg2} If $x \in Z(G)$ and $\eta(x) = 1$, then $\GK \toba(\Oc, W) = \dim W$.

\medbreak
\item\label{item:old-matiasg3} If $\dim \toba(x,W) < \infty$, necessarily $\dim W = 1$. 
\end{enumerate}
\end{lemma}

\pf \ref{item:old-matiasg1} and \ref{item:old-matiasg2} follow from Lemma \ref{lema:old-matiasg-GK}:
as $\vert G \vert$ is odd $\eta(x) \neq -1$, and $\dim W$, a divisor of $\vert G^x \vert$,
could not be 2. 
In turn \ref{item:old-matiasg1} implies \ref{item:old-matiasg3}.
\epf

\subsection{Nichols algebras with central support}\label{subsection:central-conjugacy-classes}
Let $M \in \yd{\ku G}$ with central support. 
Then the braided vector space $M$ can be realized in $\yd{\ku Z(G)}$, 
hence it fits into the theory of Nichols algebras over abelian groups sketched at the Introduction. 
For illustration we describe the semisimple $M \in \yd{\ku G}$ of finite length and central support
such that $\GK \toba(M) < \infty$, up to Conjecture \ref{conj:AAH}. See 
\cite{AAH-triang,AAH-pems} for Nichols algebras of indecomposable $M$.

\begin{prop}\label{prop:central-conjugacy-classes} Let $M \in \yd{\ku G}$ be semisimple
of the form
\begin{align}
M \simeq M_1 \oplus \cdots \oplus M_t \oplus M_{t+1} \oplus \cdots \oplus M_{\theta},
\end{align}
where $M_i \simeq M(g_i, W_i)$ with $g_i\in Z(G)$, $W_i\in \Irr G$ and 
\begin{itemize}
\item if $i\in \I_t$, then $W_i$ has dimension $\geq 2$ and central character $\eta_i$;

\medbreak
\item if $i\in \I_{t+1, \theta}$, then $W_i$ has dimension $1$ and action given by $\chi_i\in \widehat{G}$.
\end{itemize}
 Let $\bq = (q_{ij})_{i, j \in \I_{t+1, \theta}}$, where $q_{ij} = \chi_j(g_i)$.
Then $\GK \toba(M) < \infty$ if and only if the following conditions hold:

\begin{enumerate}[leftmargin=*,label=\rm{\alph*)}]
\item\label{item:central-conjugacy-classes1} The connected components of the diagram of $\bq$
 are either points labelled by 1 or else belong to the list of \cite{H-classif}.

\medbreak
\item\label{item:central-conjugacy-classes2} If $i\in \I_t$, then $\eta_i(g_i) \in \G_2 \cup \G_3$.

\medbreak
\item\label{item:central-conjugacy-classes3} If $i\in \I_t$ and $\eta_i(g_i) = 1$, 
then $\eta_i(g_j)\eta_j(g_i) = 1$ for all $i \neq j\in \I_t$  and $\eta_i(g_k)\chi_k(g_i) = 1$
for all $k \in \I_{t+1, \theta}$.

\medbreak
\item\label{item:central-conjugacy-classes4} If $i \neq j\in \I_t$ and $\eta_i(g_i) \neq  1 \neq \eta_j(g_j)$, 
then $\eta_i(g_j)\eta_j(g_i) = 1$. 

\medbreak
\item\label{item:central-conjugacy-classes5} If $i\in \I_t$, $k \in \I_{t+1, \theta}$ and $\eta_i(g_i) = \omega \in \G'_3$, 
then $\eta_i(g_k)\chi_k(g_i) = 1$ unless $\{k\}$ is a connected component of $\bq$ labelled by $-1$ and
$\eta_i(g_k)\chi_k(g_i) = \omega^2$.

\medbreak
\item\label{item:central-conjugacy-classes6} 
If $i\in \I_t$, $k \in \I_{t+1, \theta}$ and $\eta_i(g_i) = -1$, then $\eta_i(g_k)\chi_k(g_i) = 1$ except 
when $\dim W =2$ or $3$ and the points from $g_i \otimes W$ together with the connected component of $k$ 
appear in one of the following: rows 1, 8, 15 in Table 2, rows 5, 18 in Table 3, or row 8 in Table 4 from 
\cite{H-classif}.
\end{enumerate}

\end{prop}

\pf If $\GK \toba(M) < \infty$, then \ref{item:central-conjugacy-classes1} follows from \cite{H-classif}
assuming Conjecture \ref{conj:AAH}; \ref{item:central-conjugacy-classes2}, from Lemma \ref{lema:old-matiasg-GK};
\ref{item:central-conjugacy-classes3}, from Lemma \ref{lema:diagonal-ciclos} \ref{item:rossito};
\ref{item:central-conjugacy-classes4}, from Lemma \ref{lema:diagonal-ciclos} \ref{item:cycle>3}.
Now \ref{item:central-conjugacy-classes5} and \ref{item:central-conjugacy-classes6} 
follow by inspecting the list in \cite{H-classif}; the exception in \ref{item:central-conjugacy-classes5}
is from row 15, Table 2 in \cite{H-classif}. 
The proof of the converse implication is standard.
\epf

\subsection{Hopf algebras over torsion-free nilpotent groups}\label{subsection:nichols-hopf-torsion-free-main}
Combining Lemma \ref{lema:noncentral-conjugacyclasses} with Theorem \ref{thm:infinite-conjugacy-class}
and Proposition \ref{prop:central-conjugacy-classes} we get:

\begin{theorem}\label{thm:main-nichols-nilpotent-torsionfree}
Let $G$ be a finitely generated torsion-free nilpotent group.  
Let $M \in \yd{\ku G}$ be semisimple of finite length.
Then 
\begin{align*}
\GK \toba(M) \# \ku G < \infty
\end{align*}
if and only if $\supp M \subset Z(G)$ and $M$ is as
in Proposition \ref{prop:central-conjugacy-classes}. \qed
\end{theorem}

As already mentioned, the classification of all $M \in \yd{\ku G}$ of finite dimension with 
$\GK \toba(M) \# \ku G < \infty$ would follow from the abelian case once this is settled completely.

\subsection{Abelian non-central conjugacy classes}\label{subsection:noncentral-abelian-classes}

We next study Nichols algebras over finite conjugacy classes that are abelian as racks.
We keep the notation from \S \ref{subsection:first-remarks} and 
assume from now on  that 
\begin{align*}
\text{\emph{$\Oc$ is abelian but not central.}}
\end{align*}
The braided vector space $M(\Oc, W)$ can be realized in $\yd{\ku \Gamma}$ where $\Gamma = \langle \Oc \rangle$
is an abelian subgroup of $G$, hence the theory of braided vector spaces over abelian groups applies again. 
In order to give more precisions, we start with some general reductions.

\medbreak
Let $y, z \in \Oc$ and $g_y, g_z \in G$ such that $g_y \trid x = y$, $g_z \trid x = z$.
Since $\Oc$ is abelian, we see as in   \eqref{eq:thy} that 
$t_{z,y} = g_{y}^{-1}zg_y = \big(g_{y}^{-1}g_z\big) \trid x$.
Now \eqref{eq:braiding-ydG} says 
\begin{align} \label{eq:braiding-ydG-abelian-class}
&\begin{aligned}
c(g_z u\otimes g_y w) &= g_{y} (t_{z,y}\cdot w)\otimes g_z u, \\
c(g_y w\otimes g_z u) &= g_{z} (t_{y,z}\cdot u)\otimes g_y w,
\end{aligned}&  u,w &\in W.
\end{align}

We first consider the case when $\dim W \geq 2$; 
we elaborate on Lemma \ref{lema:old-matiasg-GK} using a result on pale blocks from \cite[\S 8]{AAH-triang}.

\begin{lemma}\label{lema:old-matiasg-GK-abelian-class}
Let $W \in \Irr G^x$ such that $\dim W \geq 2$ and  $\eta(x) \neq -1$. 
Then $\GK \toba(\Oc, W) < \infty$ if and only if $\dim W < \infty$, $t_{z,y}$ acts on $W$ 
by a scalar $\tau_{z,y}$ such that $\tau_{z,y} = \tau_{y,z}^{-1}$, for any $y \neq z \in \Oc$, and
either 
\begin{enumerate}[leftmargin=10pt,label=\rm{(\roman*)}]
\item\label{item:old-matiasg-abelian-class-GK1} $\eta(x) = 1$; then 
 $\toba(\Oc, W) \twist S(W^{\vert \Oc\vert})$, $\GK\toba(\Oc, W) = \dim W \vert \Oc\vert$; 

\medbreak
\item\label{item:old-matiasg-abelian-class-GK3}  or else $\eta(x) \in \G'_3$, $\dim W =2$; 
in this case $\dim\toba(\Oc, W) = 27^{\vert \Oc\vert}$.
\end{enumerate}
\end{lemma}

\pf  Fix $y \neq z \in \Oc$. Assume that $w$
is an eigenvector for $t_{z,y}$ with eigenvalue $\lambda$ and that $u \in W$ 
is an eigenvector for $t_{y, z}$ with eigenvalue $\mu$. 

\begin{enumerate}[leftmargin=*,label=\rm{(\alph*)}]
\item\label{item:old-matiasg-abelian-class-GK4}  
By \eqref{eq:braiding-ydG-abelian-class}, the 2-dimensional braided vector space spanned by $g_z u$, $g_y w$ is  of diagonal type with Dynkin diagram
$\xymatrix{ \underset{ \ }{\overset{\eta(x)}{\circ}} \ar  @{-}[r]^{\lambda \mu}
& \underset{ \ }{\overset{ \eta(x) }{\circ} }}$.

\item\label{item:old-matiasg-abelian-class-GK5}  
Assume that  $\widetilde{u} \in W$ 
satisfies $t_{y, z} \cdot \widetilde{u}=\mu(\widetilde{u} + u)$.
Let $v_1 = g_z u$, $v_2 = g_z \widetilde{u}$ and $v_3 = g_y w$.
Then the braided vector space  $V$ with basis $(v_i)_{i\in\I_3}$
has braiding given by
\begin{align*}
(c(v_i \otimes v_j))_{i,j\in \I_3} &=
\begin{pmatrix}
\eta(x) v_1 \otimes v_1&  \eta(x) v_2  \otimes v_1& \lambda v_3  \otimes v_1
\\
\eta(x) v_1 \otimes v_2 & \eta(x) v_2  \otimes v_2& \lambda v_3  \otimes v_2
\\
\mu v_1 \otimes v_3 &  \mu(v_2 +  v_1) \otimes v_3& \eta(x) v_3  \otimes v_3
\end{pmatrix}.
\end{align*}
Now \cite[Theorem 8.1.3]{AAH-triang} says that $\GK \toba(V) < \infty$ if and only if 
\begin{align}\label{eq:braiding-paleblock-point-bis}
\eta(x) &= -1 &&\text{ and } & \lambda \mu &= \begin{cases} 1, & \text{thus} \GK \toba(V)  =1; \text{ or }\\
-1, & \text{thus} \GK \toba(V)  = 2.
\end{cases}
\end{align}
\end{enumerate}

First we assume that $\GK \toba(\Oc, W) < \infty$. Since one of our assumptions is that $\eta (x) \neq -1$,
both $t_{y, z}$ and $t_{z, y}$ 
act diagonally on $W$ by \ref{item:old-matiasg-abelian-class-GK5}.

\medbreak
\noindent \ref{item:old-matiasg-abelian-class-GK1}:
If $\eta(x) = 1$, then $\dim W < \infty$ by Lemma \ref{lema:old-matiasg-GK}.  
Pick eigenvalues $\lambda$, $\mu$ as above. Then 
$\lambda \mu = 1$ by Lemma \ref{lema:diagonal-ciclos} \ref{item:rossito}; thus $\mu$ and $\lambda = \mu^{-1}$ 
are uniquely determined i.~e. $t_{z,y}$ and $t_{y,z}$ act by inverse scalars.

\medbreak
\noindent \ref{item:old-matiasg-abelian-class-GK3}:
Assume that $\eta(x) =: \omega\in \G'_3$ and $\dim W =2$.
Let $\lambda_1$ and $\lambda_2$ be the eigenvalues of $t_{y, z}$ acting on $W$, respectively $\mu_1$ and $\mu_2$
the eigenvalues of $t_{z, y}$.
Then the Dynkin diagram of $V := g_z W \oplus g_y W$ has the form
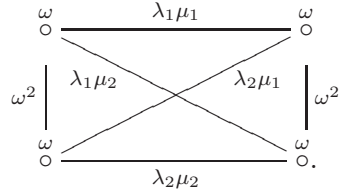
\begin{align*}
\xymatrix{ 
\underset{ \ }{\overset{\omega}{\circ}} \ar  @{-}[rrr]^{\lambda_1 \mu_1} \ar @{-}[d]_{\omega^2} 
\ar  @{-}[rrrd]_<(0.2){\lambda_1 \mu_2}
& & & \underset{ \ }{\overset{ \omega }{\circ} } \ar @{-}[d]^{\omega^2} \ar  @{-}[dlll]^<(0.2){\lambda_2 \mu_1}
\\
\underset{ \ }{\overset{\omega}{\circ}} \ar  @{-}[rrr]_{\lambda_2 \mu_2} & & & \underset{ \ }{\overset{ \omega }{\circ} }.
}
\end{align*}
If at least two of the $\lambda_i \mu_j$ are different from 1, then this diagram has either a 3-cycle  or a 4-cycle with
all vertices equal to $\omega$, so $\GK \toba(V) = \infty$, while if three of them are 1, then
the fourth also is 1. The converse  is clear.
\epf

The case when $\dim W \geq 2$ and  $\eta(x) = -1$ is still open, see  \S \ref{subsection:even-order}.

\medbreak
Next we treat the case when $\dim W = 1$, i.~e. given by a character $\chi$.

\begin{lemma}\label{lema:tot-disconnected}
Assume that $\Oc$ is an abelian rack and that
\begin{align}\label{eq:tot-disconnected}
\chi\left((g^{-1} \trid x) (g \trid x)\right) &= 1 & \text{ for every } g &\in G \backslash G^x.
\end{align}
Then $M(\Oc, \chi)$ is a braided vector space of diagonal type whose Dynkin diagram is totally disconnected
with all vertices labelled by $q = \chi(x)$. Hence
\begin{align}\label{eq:tot-disconnected-GK}
\begin{aligned}
&\dim \toba(\Oc, \chi) =  N^{\vert \Oc \vert},\  \GK\toba(\Oc, \chi) =0,
&\text{if } q &\in \G_{N}', \quad N \geq 2; \\ 
&\GK\toba(\Oc, \chi) = \vert \Oc \vert, &\text{if } q &= 1 \text{ or } q\in \ku \backslash\G_{\infty}.
\end{aligned}
\end{align}
\end{lemma}

\pf 
By \eqref{eq:braiding-ydG-abelian-class}, we have
$c(g_z \otimes g_y ) = \chi\left(\big( g_{y}^{-1}g_z\big) \trid x\right)  g_{y} \otimes g_z$, therefore 
$c^2(g_z \otimes g_y) = \chi\left(\big( g_{z}^{-1}g_y\big) \trid x\right)\chi\left(\big( g_{y}^{-1}g_z\big) \trid x\right)  g_{z} \otimes g_y$.
By \eqref{eq:tot-disconnected} the Dynkin diagram is totally disconnected and  \eqref{eq:tot-disconnected-GK}
follows.  \epf

If $\dim W = 1$ but \eqref{eq:tot-disconnected} does not hold, then we apply an argument
generalizing \cite[Lemma 3.7]{YZ}. 
Given $g \in G$, as $\Oc_x$ is finite, the set
\begin{align}\label{eq:orbg}
\orbg^x_g := \{z_i := g^{i} \trid x: i \in \Z\} \subset \Oc_x
\end{align} 
is finite;
$\vert \orbg^x_g \vert = 1$ iff $g\in G^x$ and $\vert \orbg^x_g \vert$ divides $\vert g\vert$
when this last is finite.

\begin{lemma}\label{lema:YZ}
Fix $g \in G \backslash G^x$. Assume that 
\begin{align}\label{item:YZ3}  
\zeta := \chi\left(z_{-1}z_{1}\right) = \chi\left((g^{-1} \trid x) (g \trid x)\right) \neq 1.
\end{align}

If $\GK \toba(\Oc, \chi) < \infty$, then
\begin{align}
q &=\chi\left(x\right)\neq 1,
\\ n &:= \vert \orbg^x_g\vert \in \{2,3\}.
\end{align}
Furthermore

\medbreak
\begin{enumerate}[leftmargin=*,label=\rm{(\roman*)}]
\item\label{item:YZ1} if $n = 2$, then either $\zeta = q^{-1}$ or else 
$\zeta \in \G'_{12}$ and $q = -\zeta^{2}$.

\medbreak
\item\label{item:YZ4} If $n = 3$, then  $q = -1$ and $\zeta \in \G'_3$.
\end{enumerate}
\end{lemma}

For this Lemma, we just need $\orbg^x_g$ to be an abelian subrack of $\Oc$.

\pf Let $\J = \I_{0,n-1}$.
In the notation \eqref{eq:thy}, choosing $g_{z_i} =  g^{i}$ we have that
\begin{align*}
t_{z_i,z_j} &= g_{z_j}^{-1} z_i g_{z_j} = g^{-j} z_i g^{j} = z_{i-j},&  i,j &\in \J,   
\end{align*}
since $\orbg^x_g$ is abelian.
Set $v_i = g^{i}1 \in M(\Oc, \chi)$, $V =\langle v_i: i\in \J\rangle$. Then 
\begin{align}\label{eq:braiding-zij}
c(v_i \otimes v_j) &= \chi(z_{i-j}) v_j \otimes v_i,& i,j &\in \J.   
\end{align}
Thus the Dynkin diagram of $V$ has locally the form 
$$\xymatrix{\cdots \underset{i}{\overset{q}{\circ}} \ar  @{-}[rr]^{\chi\left(z_{-1}z_{1}\right)}
& &\underset{i+1}{\overset{ q }{\circ} \cdots}}$$
By Lemma \ref{lema:diagonal-ciclos} \ref{item:rossito}, $q \neq 1$.
If $n =2$, then Remark \ref{remark:2-cycle} applies.
If $2 < n < \infty$, then this is an $n$-cycle by \eqref{item:YZ3}.
Thus  Lemma \ref{lema:diagonal-ciclos} and Remark \ref{remark:3-cycle} imply that $n=3$
and \ref{item:YZ4}.  
\epf

Assume that $G$ is finite; hence $n = \vert \orbg^x_g\vert$ divides $\vert G\vert$.
If $G$ has odd order, then \ref{item:YZ1} and \ref{item:YZ4}
could not happen and we have the following consequence.

\begin{lemma}\label{lema:YZcor}
Assume that $G$ is a finite group of odd order and that $\Oc$ is an abelian rack. 
Then $\GK \toba(\Oc, \chi) < \infty$ if and only if \eqref{eq:tot-disconnected} holds.
When this happens, $\dim \toba(\Oc, \chi) < \infty$ if and only if $q \neq 1$.
\end{lemma}

\pf If \eqref{eq:tot-disconnected} does not hold, then $\GK \toba(\Oc, \chi) = \infty$ by Lemma \ref{lema:YZ}. 
If \eqref{eq:tot-disconnected} holds, then $\GK \toba(\Oc, \chi) < \infty$ by Lemma \ref{lema:tot-disconnected}.
The last claims follows from \eqref{eq:tot-disconnected-GK} since $G$ is finite. \epf

\begin{example} If $\Oc$ is an abelian rack, then $\Oc \subset G^x$. In this case, the extra 
assumption $\Oc \subset \ker \chi$ implies
\eqref{eq:tot-disconnected} and  $q=1$. Thus $\GK \toba(\Oc, \chi) = \vert \Oc \vert$. 
\end{example}

\subsection{Nilpotent groups of odd order}\label{subsection:nichols-nilpotent-main}
We are now ready to determine the finite-dimensional Nichols algebras
over a finite nilpotent group of odd order in terms of its group structure and representation theory.

\begin{theorem}\label{thm:main-nichols-nilpotent-odd}
Let $G$ be a finite nilpotent group of odd order. Given a finite-dimensional $M \in \yd{\ku G}$,
we have that $\dim \toba(M) < \infty$ if and only if 
$M \simeq M_0 \oplus M_1 \oplus \dots \oplus M_t$ where:
\begin{enumerate}[leftmargin=10pt,label=\rm{(\roman*)}]
\medbreak
\item\label{item:main1}  $\supp M_0 \subseteq Z(G)$, 
hence $M_0$ is given by a family of YD-pairs $(g_i, \chi_i)_{i\in J}$
such that the connected components of the matrix $\bq = (q_{ij})_{i,j \in J}$ 
belong to the list in \cite{H-classif}.

\medbreak
\item\label{item:main2} For $j\in \I_t$, $M_j \simeq M(\Oc_j, \chi_j)$ where $\Oc_j$ is not central 
and  abelian as rack; 
$\chi_j \in \widehat{G^{x_j}}$ for a fixed $x_j \in \Oc_j$ that satisfies \eqref{eq:tot-disconnected}; 
and $q_j :=\chi_j\left(x_j\right)$ has order $2 < N_j < \infty$.
Also $\dim \toba(\Oc_j, \chi_j) = N_j^{\vert \Oc_j \vert}$.
\end{enumerate}
Furthermore,
\begin{align}\label{eq:braiding-square-trivial}
c^2_{\vert M_i \otimes M_j} &= \id_{\vert M_i \otimes M_j},& i\neq j&\in \I_{0,t}.
\end{align}

\end{theorem}

\pf Since $\yd{\ku G}$ is semisimple, we may decompose $M = \oplus_{i \in I} M_i$ where the $M_i$'s are simple.
Let $J = \{i\in I: \supp M_i \subseteq Z(G)\}$ and $M_0 = \oplus_{i\in J} M_i$.
If $\supp M_i$ is central, then $M_i \simeq M(\{g_i\}, W)$ where $g_i\in Z(G)$;
by Lemma \ref{lema:old-matiasg} \ref{item:old-matiasg3}, $\dim W = 1$.
Thus there is a family of YD-pairs $(g_i, \chi_i)$ such that
$M_i \simeq \ku_{g_i}^{\chi_i}$ for all $i \in J$. By \cite{H-classif}, 
$\dim \toba(M_0) < \infty$ if and only if \ref{item:main1} holds.

\medbreak
Assume that $M \simeq M(\Oc, W)$ for some non-central conjugacy class $\Oc$. 
By Theorem \ref{thm:odd-order-abelian-or-typeC} $\Oc$ is either of type C, or  is an abelian rack.
In the first case $\dim \toba(M) = \infty$ by Theorem \ref{th:typeC}. 
Assume that $\Oc$ is abelian and $\dim \toba(M) < \infty$; fix $x \in \Oc$.
By Lemma \ref{lema:old-matiasg}, $\dim W = 1$, so $W$ is given by $\chi \in \widehat{G^x}$.
By Lemma \ref{lema:YZcor}, \eqref{eq:tot-disconnected} should hold. 
Conversely, if $\Oc$ is abelian, $W$ is given by $\chi \in \widehat{G^x}$ and \eqref{eq:tot-disconnected} 
holds, then $\dim \toba(M) < \infty$ by Lemma \ref{lema:YZcor}.
Up to renumbering $I \backslash J$, we have \ref{item:main2}.

\medbreak
Finally, assume that $M', M'' \in \Irr \yd{\ku G}$ satisfy $\dim \toba(M'\oplus M'') < \infty$, 
where $M' \simeq M(\Oc, \chi)$ and $\Oc$ is a non-central conjugacy class.
By \ref{item:main2}, $M'$ is of diagonal type.
We may replace $G$ by the subgroup generated by $\supp M'\cup \supp M''$.
Clearly $\dim M'$ is odd and $>1$. If \eqref{eq:braiding-square-trivial} does not hold, then
$\dim M'$ should be 3 by Theorem \ref{th:HV}, but then $M'$ is not of diagonal type, a contradiction.
\epf

\begin{algorithm}\label{algor:nichols-nilpotent-odd} Let $G$ be a finite nilpotent group of odd order. 
By Theorem \ref{thm:main-nichols-nilpotent-odd}, to list all Nichols algebras in $\yd{\ku G}$
 we should do the following.

\begin{enumerate}[leftmargin=10pt,label=\rm{(\roman*)}]
\medbreak
\item\label{item:algorithm1}  
Compute $Z(G)$, $[G,G]$ and $\widehat{G} = \widehat{[G, G]}$. Thus we have 
all YD-pairs $(g, \chi)\in Z(G)\times \widehat{G}$. 

\medbreak
\item\label{item:algorithm1.5}
For any braided vector space $(V, c^{\bq})$ either of Cartan type or of diagonal type \eqref{eq:dynkin-br(2,a)}
or else \eqref{eq:dynkin-br(3)}, compute all principal realizations over $G$, see
Definition \ref{remark:ppal-realiz-odd-order} and the subsequent discussion.

\medbreak
\item\label{item:algorithm2}  
Compute the set of abelian conjugacy classes
\begin{align}
\Cla (G) =\{\Oc \text{ conjugacy class of }G: [\Oc, \Oc]  =e, \ \Oc \not\subset Z(G) \}.
\end{align}
Given $\Oc \in \Cla(G)$, pick $x \in \Oc$ and compute $G^x$ and $\widehat{G^x}$. Thus
the set of pairs (with an evident abuse of notation)
\begin{align}\label{eq:pairs-YDmodules}
\{(\Oc, \chi): \Oc \in \Cla (G), \chi \in \widehat{G^x} \text{ satisfies \eqref{eq:tot-disconnected} and } \chi(x) \neq 1 \}
\end{align}
parametrizes the finite-dimensional Nichols algebras in $\yd{\ku G}$ of irreduci\-ble objects
with abelian non-central support.

\medbreak
\item\label{item:algorithm3}
Compute all pairs $(\Oc, \chi), (\Oc', \chi')$ as in \eqref{eq:pairs-YDmodules} such that 
\begin{align*} 
[\Oc, \Oc']  &= e,
\\ \chi((g'_z)^{-1} \trid y) \chi'((g_y)^{-1} \trid z) &= 1,& \forall y&\in \Oc, z \in \Oc',
\end{align*}
where $x\in \Oc$, $x' \in \Oc'$ and $g_y, g'_z \in G$ satisfy $g_y \trid x = y$, $g'_z \trid x' = z$.

\medbreak
\item\label{item:algorithm4}
Similarly for $(\Oc, \chi)$ as in \eqref{eq:pairs-YDmodules} and $(g', \chi')\in Z(G)\times \widehat{G}$. 

\end{enumerate}
\end{algorithm}

\subsection{Nilpotent groups of odd order, II}\label{subsection:nichols-nilpotent-finiteGK}
Let $G$ be a finite nilpotent group of odd order. We extend the discussion in the previous subsection
to determine all $M\in \yd{\ku G}$ such that $\GK \toba(M) < \infty$.  
To have a complete picture we still need Conjectures \ref{conj:AAH} and \ref{conjecture:typeC-X-finite-GK}, and 
an analogue of Theorem \ref{th:HV}. Here are the necessary steps:

\begin{algo} $\supp M \subset Z(G)$. \end{algo} 
Since every object in $\yd{\ku G}$ is semisimple, 
Proposition \ref{prop:central-conjugacy-classes} gives a complete picture, up to Conjecture \ref{conj:AAH}.

\begin{algo}
$M \simeq M(\Oc, W)$ where $\Oc$ is not central.
\end{algo}
 By Theorem \ref{thm:odd-order-abelian-or-typeC}, $\Oc$ is either abelian or of type C.
To discard type C, we need the validity of Conjecture \ref{conjecture:typeC-X-finite-GK}.

Assume that $\Oc$ is abelian; fix $x \in \Oc$. If $\dim W > 1$, then $\eta(x) = 1$ and $\GK \toba(\Oc, W) = \dim W \vert \Oc\vert$ by Lemma \ref{lema:old-matiasg-GK-abelian-class}. 
If $\dim W = 1$, then \eqref{eq:tot-disconnected} should hold and then $\GK \toba(M) < \infty$, given by \ref{eq:tot-disconnected}, see Lemma \ref{lema:YZcor}.

\begin{algo} Braidings between $M(\Oc, W)$ with $\Oc$ not central and other summands.
\end{algo}
We guess that \eqref{eq:braiding-square-trivial} holds; this would 
need an analogue of Theorem \ref{th:HV} for finite $\GK$. Clearly this would be related to Conjecture \ref{conjecture:typeC-X-finite-GK}.

\subsection{Finitely generated nilpotent groups whose torsion has order coprime to 6}\label{subsection:nichols-nilpotent-finiteGK-torsion-coprime-6}
Let $G$ be a finitely generated nilpotent group with torsion subgroup $T$,
let $\Oc$ be a finite abelian conjugacy class and $x \in \Oc$.
Assume that $g \in G \backslash G^x$ satisfies \eqref{item:YZ3}; 
recall the set $\orbg^x_g$ from \eqref{eq:orbg}.

\begin{lemma}\label{lema:orbg-divides-torsion}
If $p = \vert \orbg^x_g \vert$ is prime, then it divides $\vert T\vert$.
\end{lemma}

\pf Let $K = \langle g, G^x \rangle$. As $g^p \in G^x$,
$[K: G^x]\in p^{\N}$ by Theorem \ref{theorem:malcev-finite-index}, thus 
$p$ divides $\vert \Oc\vert = [G: G^x] = [G: K][K: G^x]$.
Then Lemma \ref{lema:order-conjugacy-class} applies.
\epf

Assume now that $\vert T\vert$ is coprime to 6.  
In particular $\vert T\vert$ is odd, hence any finite conjugacy class is either abelian or of type C by
Proposition \ref{prop:noncentral-conjugacyclasses-torsion-typeC}.

\begin{lemma}\label{lema:nichols-nilpotent-coprime-6-dimW1}
Let  $\chi \in \widehat{G^x}$. 
Then  $\GK\toba(\Oc, \chi) < \infty$ if and only if \eqref{eq:tot-disconnected} holds, in which case 
$\toba(\Oc, \chi)$ is given by \eqref{eq:tot-disconnected-GK}.
\end{lemma}

\pf If \eqref{eq:tot-disconnected} does not hold, then $\GK\toba(\Oc, \chi) = \infty$ by Lemmas 
\ref{lema:YZ} and \ref{lema:orbg-divides-torsion} (that excludes $2$ and $3$ by the hypothesis).
The converse is clear.
\epf

In summary, if $G$ is a finitely generated nilpotent group with torsion $T$ with $\vert T\vert$ coprime to 6,
and assuming Conjecture \ref{conjecture:typeC-X-finite-GK}, 
then the $M = M(\Oc, W) \in \Irr \yd{G}$  such that $\GK \toba(M) < \infty$ are 
either covered by Proposition \ref{prop:central-conjugacy-classes}, 
Lemma \ref{lema:old-matiasg-GK-abelian-class} and Lemma \ref{lema:nichols-nilpotent-coprime-6-dimW1}, or else
$\Oc$ is abelian non-central, $\dim W \geq 2$ and $\eta(x) = -1$.
For these last, see the discussion in \S \ref{subsection:even-order}.
Once this is settled, the determination of the semisimple $M$ such that $\GK \toba(M) < \infty$
can be obtained softly. 

\bigbreak
\subsection{Examples}\label{section:examples}

\subsubsection{Class 2}  
Let $G$ be a finite nilpotent group of odd order. We discuss how to prove the following result 
from \cite{YZ} by means of Theorem \ref{thm:main-nichols-nilpotent-odd}: \emph{If
$[G, G] = Z(G)$, then $\dim\toba(M) = \infty$ for any $M \in \yd{\ku G} \backslash 0$}.
Equivalently any finite-dimensional pointed Hopf algebra over $G$ is isomorphic to $\ku G$.

\smallbreak
The contention $[G, G] \supset Z(G)$ implies that $\dim\toba(M) = \infty$ for any $M \in \yd{\ku G}$
with central support. By Theorem \ref{thm:main-nichols-nilpotent-odd}, we are reduced to prove:

\setcounter{claim}{0}
\renewcommand{\theclaim}{\arabic{claim}}
\begin{claimo}\label{claim:nilpotent-notcentral-1}
Let $x\in G$ such that $\Oc_x$ is non-central abelian and let
$\chi \in \widehat{G^{x}}$ with $\chi\left(x\right) \neq 1$.
Then there exists $g \in G \backslash G^x$ such that $\chi\left((g^{-1} \trid x) (g \trid x)\right) \neq 1$.
\end{claimo}
First, let $\Gamma$ be a group, 
$x \in \Gamma$ and $g \in N_\Gamma(\Gamma^x)$, such that $[g, [g, x]] = 1$.
Then  $(g^{-1} \trid x) (g \trid x) = x^2$.
Indeed, clearly $[g, x]^{-1} = [g^{-1}, x]$, hence
\begin{align*}
(g^{-1} \trid x) (g \trid x) = [g^{-1}, x]x [g, x]x = [g, x]^{-1} [g, x]x^2 = x^2.
\end{align*}

Now $N_G(G^x) \neq G^x$ since $G$ is nilpotent and any
$g \in N_G(G^x) \backslash G^x$  satisfies $[g, [g, x]] = 1$ because $[G, G] \subset Z(G)$.

\subsubsection{Heisenberger groups}\label{subsection:nichols-heisenberger3}
Let $K$ be a commutative ring and $n \in \N$.  We consider
the Heisenberg group $\mathtt{H} = \heis{2n+1}{K}$ that consists of the matrices
\begin{align}\label{eq:heis-matrix}
 \begin{pmatrix}
1 & a_{1}& a_{2} &\dots &a_{n}  & c\\
0 & 1& 0 &\dots & 0 &  b_{1}\\
0 & 0& 1 &\dots & 0& b_{2}\\
\vdots & \vdots & \vdots & \ddots   & \vdots \\
0 & 0& 0 &\dots & 1& b_{n} \\
0 & 0 & 0 &\dots & 0& 1
\end{pmatrix} \in GL_{n+2}(K).
\end{align}
For simplicity, we denote $\ba := (a_{1}, a_{2}, \dots, a_{n})$, 
$\bb := (b_{1}, b_{2}, \dots, b_{n})$ and the matrix \eqref{eq:heis-matrix} by 
$(\ba, \bb, c)$.
Let $\omega: K^{2n} \times K^{2n} \to K$ be the `symplectic form' 
\begin{align*}
\omega\left((\ba, \bb) , (\br, \bs)\right) & =  \sum_{1 \le i \le n}(a_is_i - b_ir_i),&
\ba, \bb, \br, \bs &\in K^{n}.
\end{align*}
Then
$(\ba, \bb, c) \trid (\br, \bs, t)  = \big(\br, \bs,  t + \omega\left((\ba, \bb) , (\br, \bs)\right)\big)$.

Given $\br, \bs \in K^{n}$, let $(\br, \bs)^{\perp}$ be the $K$-submodule of $K^{2n}$ 
of those $(\ba, \bb)$ such that
$\omega\left((\ba, \bb) , (\br, \bs)\right) = 0$; and let
$\langle \br, \bs\rangle$ be the 
ideal of $K$ generated by $r_{i}, s_{i}$, $i \in \I_{n}$. 
Clearly $Z(\mathtt{H}) = 0\times 0\times K$ and  the conjugacy class of $(\br, \bs, t)$ is
\begin{align*}
\Oc_{(\br, \bs, t)} &:=\{(\br, \bs, t + \ell ): \ell \in \langle \br, \bs\rangle \}.
\end{align*} 
Fix a non-central class $\Oc$ and $(\br, \bs, t) \in \Oc$, where 
$(\br, \bs) \neq (0,0)$. Given  $\ell \in \langle \br, \bs\rangle$
pick $(\ba, \bb)$ such that $\omega\left((\ba, \bb) , (\br, \bs)\right) = \ell$ 
and set 
\begin{align*}
x_\ell &= (\br, \bs, t+\ell) \in \Oc, & g_\ell &= (\ba, \bb, 0),
\end{align*}
so that $g_\ell \trid x_0 = x_\ell$. Then for $m\in\langle \br, \bs\rangle$, the element \eqref{eq:thy} is given by
\begin{align*}
t_{m,\ell} := g_\ell^{-1}x_m g_\ell &= (-\ba, -\bb, \sum_{i}a_ib_i) (\br, \bs, t+m) (\ba, \bb, 0) = (\br, \bs, t+ m-\ell).
\end{align*}

The centralizer of $(\br, \bs, t)$ is $\mathtt{H}^{(\br, \bs, t)} \simeq (\br, \bs)^{\perp}\times K$, 
which is abelian, so  $\widehat{\mathtt{H}^{(\br, \bs, t)}} \simeq \widehat{(\br, \bs)^{\perp}}\times \widehat{K}$ parametrizes $\Irr \mathtt{H}^{(\br, \bs, t)}$.  Let $\chi = (\chi_1, \chi_2)$ be such a character and $q = \chi_1(\br, \bs)\chi_2(t)$. 
Then the braiding of $M(\Oc, \chi)$ is given by 
\begin{align}\label{eq:braiding-ydHeis}
c(g_m \otimes g_\ell) &= \chi_1(\br, \bs) \chi_2(t) \chi_2(m-\ell) g_{\ell}  \otimes g_m,& m,\ell &\in \langle \br, \bs\rangle.
\end{align}
That is,  $M(\Oc, \chi)$ is  of diagonal type with matrix twist-equivalent to 
$\bq = (q_{ij})_{i, j \in \I}$,  $q_{ij} = q$ for all $i,j$. Clearly $\dim M(\Oc, \chi) = \vert \langle \br, \bs\rangle \vert$.
By Lemma \ref{lema:diagonal-q-constante}, $\GK \toba (\Oc, \chi) =$
\begin{align*}
  &= \begin{cases}
 \vert \langle \br, \bs\rangle \vert, & \toba (\Oc, \chi)\twist S(M(\Oc, \chi)) \text{ if }
\chi_1(\br, \bs) = \chi_2(t)^{-1}; \\
0, &\toba (\Oc, \chi) \twist \Lambda(M(\Oc, \chi)) \text{ if }
\chi_1(\br, \bs) = -\chi_2(t)^{-1}; \\
0, & \dim\toba(x, \chi) = 27 \quad
\text{if }\chi_1(\br, \bs) \chi_2(t) \in \G'_3 \text{ and } \vert \langle \br, \bs\rangle \vert = 2;
 \\
\infty &\text{in any other case.}
\end{cases}
\end{align*}
The first two lines correspond to the situation in Lemma \ref{lema:tot-disconnected}
while the  third one is covered by Lemma \ref{lema:YZ} \ref{item:YZ1}; this last case  occurs e.~g.
with $K = \Z/6$, $n =1$, $\br = (3,0) = \bs$, $t = \overline 2$, $\chi_1$ trivial and $\chi_2$ non-trivial.  

\bigbreak
Another family of examples arises taking an ideal $I$ of $K$ and 
the quotient $\widetilde{\mathtt{H}} := \heis{2n+1}{K}/ 0\times 0\times I \simeq K^{2n} \times K/I$.
Clearly $Z(\widetilde{\mathtt{H}})  =  I^{2n} \times K/I$.
Let $\pi:K \to K/I$,
$t \mapsto \overline t$ be the natural projection and
let $\overline \omega: K^{2n} \times K^{2n} \to K/I$, 
$\overline \omega = \pi\omega$. If $(\br, \bs) \in K^{2n}$,  and $t \in K$, then
\begin{align*}
\widetilde{\mathtt{H}}^{(\br, \bs, \overline t)}&= (\br, \bs)^{\perp,I} \times K/I;
&
\Oc_{(\br, \bs, \overline t)} &= (\br, \bs, \overline{ t + \langle \br, \bs\rangle}),
\end{align*}
where $(\br, \bs)^{\perp,I} = \{(\ba, \bb ) \in K^{2n}: \, \omega\left((\ba, \bb) , (\br, \bs)\right) \in I\}$.
Let $\chi = (\chi_1, \chi_2) \in \widehat{(\br, \bs)^{\perp, I}}\times \widehat{K/I}$  and $q = \chi_1(\br, \bs)\chi_2(t)$.
Then $M(\Oc, \chi)$ is  of diagonal type with matrix twist-equivalent to 
$\bq = (q_{ij})_{i, j \in \I}$,  $q_{ij} = q$ for all $i,j$ and $\toba(\Oc, \chi)$ is determined by a similar analysis.

Suppose that $K = \Z$ and $I = \Z/N$ where $N \geq 2$; then $\mathtt{H}$ is an FC-group with torsion $\simeq \Z/N$.
Assume further that $N = 2d$ is even, $n =1$, $\br = (d,0) \neq \bs = (0,0)$, $t =  \overline 1$; thus
$(\br, \bs)^{\perp, I} \simeq \frac{N}{(d, N)}\Z \times \Z^3 \simeq \Z^4$.
Then $\vert \overline{\langle \br, \bs\rangle}\vert = 2$ and $q = \chi_1(\br, \bs)\chi_2(t)$ might be in $\G'_3$
even if $N$ is coprime to 3.

\subsubsection{Unitriangular groups}\label{section:nichols-triangular3}
We show examples of classes of type C.
Let $K$ be a commutative ring and ${\mathtt G} = \ut{4}{K}$, consisting  of the matrices
\begin{align}\label{eq:ut-matrix}
\begin{pmatrix}
1 & a_{12} & a_{13} & a_{14}\\
0 & 1 & a_{23} & a_{24} \\
0 & 0 & 1 & a_{34}\\
0 & 0 & 0 & 1
\end{pmatrix} \in GL_{4}(K).
\end{align}
For simplicity, we denote matrices \eqref{eq:ut-matrix} by 
$\ba := \begin{pmatrix}a_{12} \qquad a_{23} \qquad  a_{34}\\
a_{13}   \qquad a_{24}  \\
a_{14} \end{pmatrix}$, etc. Then
\begin{align*}
\ba \trid \bb & =
\begin{pmatrix} b_{12} \qquad\qquad\qquad b_{23} \qquad\qquad\qquad  b_{34}\\
 b_{13} + a_{12}b_{23} - b_{12}a_{23}   \qquad\qquad b_{24} + a_{23}b_{34} - b_{23}a_{34}  \\
b_{14} + b_{12}(a_{23}a_{34} - a_{24}) - b_{13}a_{34} - b_{23}a_{12}a_{34} +b_{24}a_{12} + b_{34}a_{13} \end{pmatrix}.
\end{align*}

Let $\br := \begin{pmatrix}r_{12} \qquad r_{23} \qquad  r_{34}\\
0\qquad 0  \\ 0 \end{pmatrix}$
where $r_{12}, r_{23}, r_{34} \in \Uc(K)$ (the group of units of $K$). Then 
\begin{align*}
\Oc_{\br} = \left\{ \begin{pmatrix}r_{12} \qquad r_{23} \qquad  r_{34}\\
c_{13}\qquad c_{24}  \\ c_{14} \end{pmatrix}:
c_{13}, c_{24}, c_{14} \in K
\right\}.
\end{align*}
Then $\Oc_{\br}$  is not an abelian, 
e.~g. $\bs := \begin{pmatrix}r_{12} \qquad r_{23} \qquad  r_{34}\\
1\qquad 0  \\ 0 \end{pmatrix}$ satisfies
$[\br, \bs] = \begin{pmatrix}
0 \qquad  0 \qquad 0 \\
0 \qquad 0  \\
- r_{34} 
\end{pmatrix}$.
Hence $H := \langle \br, \bs \rangle \simeq \heis{3}{K}$ via $(1,0,0) \mapsto \br$, 
$(0,1,0) \mapsto \bs$, and $\Oc_{\br}^H \cap \Oc_{\bs}^H = \emptyset$.
Thus, if $K$ is finite and $\vert K \vert$ is odd, then $\Oc_{\br}$ is of type C.

\section{Conclusions}\label{section:conclusions}
We discuss the problems that remain open as well as some applications.

\subsection{Representations}\label{subsection:nichols-nilpotent-representations}
In the previous sections we have argued assuming some information about the representations and conjugacy classes
of finitely generated nilpotent groups.
But the representation theory of such groups is not completely known to our knowledge. 
Let $H \leq G$ and $W\in \Rep H$. 
Then $\Ind_{H}^{G} W = \ku G \otimes_{\ku H} W$ is called the induced representation of $\rho$.
We need two definitions. Let $\pi: G \to V$ be a representation.

\smallbreak
\noindent $\circ$ $\pi$ is \emph{monomial} if there 
are $H \leq G$ and $\chi \in  \widehat{H}$ such that $\pi \simeq \Ind_{H}^{G} \chi$.

\smallbreak
\noindent $\circ$ $\pi$  \emph{has finite weight} if there are  $K \leq G$ and $\eta \in  \widehat{K}$ such that 
\begin{align*}
V_{K, \eta} := \{v \in V: \pi(k)(v) = \eta(k)v \ \forall k \in K \}
\end{align*}
is non-zero and finite-dimensional.
Generalizing classical results of Dixmier, Kirillov and 
Brown, Parshin \cite{parshin} conjectured the following result.

\begin{theorem}\cite{beloshapka-gorchinskiy}.
Let $G$ be a finitely generated nilpotent group.
An irreducible representation of $G$ is monomial iff it has a finite weight.
\end{theorem}

Thus finite-dimensional irreducible representations of $G$ are monomial, 
but this was already known \cite[Lemma 1]{brown}.

\subsection{Finite nilpotent groups of even order}\label{subsection:even-order}
In order to deal with Nichols algebras
over a finite nilpotent group of even order   with finite dimension or finite $\GK$  
extending Theorem \ref{thm:main-nichols-nilpotent-odd} (see also Section \ref{subsection:nichols-nilpotent-finiteGK}),    
the following points need to be addressed:

\begin{itemize}[leftmargin=10pt,label=$\circ$]
\item \emph{Conjugacy classes}. It suffices to consider $2$-groups.
We need to keep track of the conjugacy classes that are neither abelian nor of type C;
the information from \cite{HV,HV-rank>2} would be crucial. 

\medbreak
\item \emph{Irreducible Yetter-Drinfeld modules, $\dim W \geq 2$}. 
In the setting of Lemma \ref{lema:old-matiasg-GK-abelian-class} we still have to consider the case $\eta(x) = -1$.
We shall need:
\end{itemize}
\begin{lemma}
Let $\Gamma$ be an abelian group and $V = V_g \oplus V_h \in \yd{\ku \Gamma}$ such that $\GK \toba(V) < \infty$,
where $g \neq h$.
Then the action  of $h$ on $V_g$ is locally finite. \qed
\end{lemma}

Therefore we may consider eigenvectors and blocks, see the proof of Lemma \ref{lema:old-matiasg-GK-abelian-class}.
Since the classification of the Nichols algebras of \emph{blocks \& pale blocks \& points} with finite $\GK$
is not yet finished, we could not carry out a complete analysis that eventually would be straightforward.  

 \begin{itemize}[leftmargin=10pt,label=$\circ$]
\item \emph{Ditto, $\dim W = 1$}. 
We need to explore further the possibilities $\vert \orbg^x_g\vert \in \{2,3\}$ in Lemma \ref{lema:YZ},
together with other restrictions.

\medbreak
\item \emph{Braiding between simple modules}.  
The condition \eqref{eq:braiding-square-trivial} has to be adapted, cf. Theorem \ref{th:HV}; in the even case 
it does not follow from finite dimension.

\end{itemize}

\subsection{Finitely generated nilpotent groups with even torsion}\label{subsection:nilpotent-even-torsion}
The strategy would be parallel to the one in the previous Subsection.

\subsection{Finitely generated nilpotent-by-finite groups}\label{subsection:nilpotent-by-finite}

Assume that $G$ is nilpotent-by-finite and fix
$N \lhd G$ a normal nilpotent subgroup of finite index. Then the image $\widetilde{\Oc}$ of $\Oc$ in $G/N$
is a conjugacy class. As said, the knowledge of the Nichols algebras over finite groups
is still incomplete, but if $\widetilde{\Oc}$ is of type C, D or F, then so would be $\Oc$
and then any Nichols algebra over $G$ with support $\Oc$ would have infinite dimension, 
and infinite $\GK$, if  Conjecture \ref{conjecture:typeC-X-finite-GK} and its analogues for types D and F are true. See \cite{ACG-III} and references therein.

\subsection{Applications}\label{subsection:nichols-nilpotent-applications}

Once all the Nichols algebras over $G$ with finite $\GK$ or finite dimension are known, 
one still needs to (i) compute all post-Nichols algebras with finite $\GK$  and (ii)
compute all liftings. For (i), when the braided vector spaces come from the abelian setting, 
we know: a finite-dimensional Nichols algebra
does not have post-Nichols algebras with finite dimension (except itself) \cite{angiono-diagonal}. For finite $\GK$
see \cite{A Sanmarco,angiono-campagnolo-sanmarco}. As for (ii), see \cite{angiono-garcia}
Finally observe that when $G$ is torsion free, our results contribute to the classification
pointed Hopf algebras with finite $\GK$ that are domains.


\begin{thebibliography}{99}
\bibitem{A-leyva} N. Andruskiewitsch. \emph{An Introduction to Nichols Algebras}. 
In Quantization, Geometry and Noncommutative Structures in Mathematics and Physics. 
A. Cardona et al., eds., pp. 135--195, Springer (2017).


\bibitem{AA-diag} N. Andruskiewitsch,  I. Angiono. 
\emph{On finite-dimensional Nichols algebras of diagonal type}. 
Bull. Math. Sci. {\bf 7}, 353--573 (2017).

\bibitem{AAH-triang} 
N. Andruskiewitsch,  I. Angiono,  I. Heckenberger.
\emph{On finite GK-dimen\-sional Nichols algebras over abelian groups}. Mem. Amer. Math. Soc., to appear.

\bibitem{AAH-diag} \bysame
\emph{On finite GK-dimensional Nichols algebras of diagonal type}. 
Contemp. Math. \textbf{728} 1--23 (2019).


\bibitem{AAH-infinite-rank} \bysame
\emph{On Nichols algebras of infinite rank with finite Gelfand–Kirillov dimension}. 
Rend. Lincei Mat. Appl. \textbf{31} 81--101 (2020).

\bibitem{AAH-pems} \bysame
\emph{Liftings of Jordan and super Jordan planes}. Proc. Edinb. Math. Soc., II. Ser. \textbf{6} 661--672 (2018).  


\bibitem{aam} N. Andruskiewitsch,  I. Angiono,  M. Moya Giusti. \emph{Nichols algebras of pale blocks}. 
In preparation.

\bibitem{ACG-III}  N. Andruskiewitsch, G. Carnovale, G. A. Garc\'ia. 
\emph{Finite-dimensional pointed Hopf algebras over finite simple groups of Lie type  III. 
Semisimple classes in $PSL_n(q)$}, Rev. Mat. Iberoam.  \textbf{33(3)}, 995--1024, (2017).




\bibitem{AG-adv} N. Andruskiewitsch and M. Gra\~na.
\emph{From racks to pointed Hopf algebras},
Adv. Math.  \textbf{178}  (2003), 177--243.



\bibitem{A Sanmarco} N. Andruskiewitsch, G. Sanmarco.
\emph{Finite GK-dimensional pre-Nichols algebras of quantum linear spaces and of Cartan type}. 
Trans. Am. Math. Soc., to appear.


\bibitem{AS Pointed HA} N. Andruskiewitsch, H.-J. Schneider. \emph{Pointed Hopf algebras}, 
New directions in Hopf algebras, MSRI series, Cambridge Univ. Press; 1--68 (2002).

\bibitem{AS-crelle}
\bysame \emph{A characterization of quantum groups}. 
J. Reine Angew. Math. {\bf 577} (2004), 81--104.

\bibitem{angiono-convex} I. Angiono.
\emph{A presentation by generators and relations of Nichols algebras of diagonal type and convex orders on root systems}.
J. Europ. Math. Soc. {\bf 17}, 2643--2671 (2015).

\bibitem{angiono-diagonal} \bysame
\emph{On Nichols algebras of diagonal type}. 
J. Reine Angew. Math. {\bf 683}, 189--251 (2013).

\bibitem{angiono-campagnolo-sanmarco} I. Angiono, E. Campagnolo, G. Sanmarco.
\emph{Finite GK-dimensional pre-Nichols algebras of super and standard type}.
\texttt{arxiv.org/2009.04863}

\bibitem{angiono-garcia} I. Angiono,  A. Garc\'ia Iglesias.
\emph{Liftings of Nichols algebras of diagonal type II. All liftings are cocycle deformations}.
Sel. Math. New Ser. 25, 5 (2019).

\bibitem{angiono-garcia2} \bysame
\emph{On finite GK-dimensional Nichols algebras of diagonal type: rank 3}. 
to appear.

\bibitem{baer}
R. Baer. \emph{Finiteness properties of groups}. Duke Math. J., \textbf{15}, 1021--1032 (1948).


\bibitem{beloshapka-gorchinskiy}
 I. V. Beloshapka, S. O. Gorchinskiy, 
 \emph{Irreducible representations of finitely generated nilpotent groups}. Mat. Sb. \textbf{207:1}  45--72 (2016); 
 Sb. Math., \textbf{207:1} 41--64 (2016).


\bibitem{brown} I. D. Brown.
{\em Representation of Finitely Generated Nilpotent Groups}.
Pac. J. Math. \textbf{45}, 13--26 (1973).

\bibitem{brown-couto}
K. A. Brown and M. Couto.
\emph{Affine commutative-by-finite Hopf algebras}. J. Algebra \textbf{573} 56--94 (2021).


\bibitem{brown-zhang}
K. A. Brown and J. J. Zhang.
 \emph{Survey on Hopf algebras of GK-dimension 1 and 2}. Contemp. Math., to appear.





\bibitem{Carnovale-Costantini-VI}  G. Carnovale, M. Costantini.
\emph{Finite-dimensional pointed Hopf algebras over finite simple groups of Lie type  VI. Suzuki and Ree groups}, 
J. Pure Appl. Alg. \textbf{225} (2021).

\bibitem{CMZ} 
A. E. Clement, S. Majewicz, M. Zyman.
\emph{The theory of nilpotent groups}. Cham: Birkh\"auser/Springer.
xvii + 307, (2017).

\bibitem{eckhardt-gillaspy} 
C. Eckhardt, E. Gillaspy. 
\emph{Irreducible representations of nilpotent groups
generate classifiable $C^\ast$-algebras}
M\"unster J. Math. \textbf{9} 253--261 (2016).


\bibitem{Grana}  M. Gra\~na.
\emph{On Nichols algebras of low dimension}. Contemp. Math. \textbf{267}, 111--134 (2000).

\bibitem{gromov} M. Gromov.
\emph{Groups of polynomial growth and expanding maps. Appendix by Jacques Tits}. 
Publ. Math., Inst. Hautes Étud. Sci. \textbf{53}, 53--78 (1981).
 
\bibitem{gruenberg} K. W. Gruenberg. \emph{Residual properties of infinite soluble groups}. 
Proc. Lond. Math. Soc. (3) 7, 29--62 (1957).

\bibitem{H-classif}  I. Heckenberger.  {\em Classification of arithmetic
root systems}, Adv. Math. \textbf{220} (2009), 59--124.

\bibitem{HY-groupoid} I. Heckenberger, H. Yamane. 
\emph{A generalization of Coxeter groups, root systems and Matsumoto's theorem}.
Math. Z. \textbf{259}, 255--276 (2008).


\bibitem{HV}
I. Heckenberger, L. Vendramin. 
\textit{The classification of Nichols algebras with finite root system of rank two}.
J. Europ. Math. Soc. \textbf{19}, 1977--2017 (2017).

\bibitem{HV-rank>2} \bysame
\emph{A classification of Nichols algebras of semi-simple Yetter-Drinfeld modules over non-abelian groups}.
J. Europ. Math. Soc. \textbf{19},  299--356 (2017).

\bibitem{KL} G. Krause and T. Lenagan,  Growth of algebras and Gelfand-Kirillov dimension. Revised ed.
Graduate Studies in Math. \textbf{22}. Amer. Math. Soc. (2000).


\bibitem{parshin}
 A. N. Parshin. \emph{Representations of higher adelic groups and arithmetic}. In
Proc. of the ICM 2010, Hyderabad.
Vol. 1. (Hindustan Book Ag., New Delhi, 2010), pp. 362--392.


\bibitem{Rad-libro} Radford, D. E., \emph{Hopf algebras}. World Scientific. xxii, 559 p.  (2012).



\bibitem{robinson} D. J. S. Robinson. 
\emph{A course in the theory of groups}. 2nd ed.  Springer-Verlag (1995).



\bibitem{R quantum groups}  M. Rosso. 
\emph{Quantum groups and quantum shuffles}. 
Invent. Math. \textbf{133}, 399--416 (1998). 





\bibitem{tomkinson} M. J. Tomkinson.
\emph{FC-groups}.
Research Notes in Math. 96. Pitman  (1984).


\bibitem{wallach}  N. Wallach. 
\emph{Real Reductive Groups I}. Academic Press (1988).



\bibitem{YZ} Y. Yang, S. Zhu. \emph{Finite-dimensional Nichols algebras over certain p-groups}. 
Commun. Alg. 45(9), 3691--3702 (2016).

\end{thebibliography}
\end{document}